\documentclass{amsart}
\usepackage{amsfonts}

\newtheorem{theorem}{Theorem}
\theoremstyle{plain}

\newtheorem{corollary}{Corollary}

\newtheorem{definition}{Definition}
\newtheorem{example}{Example}

\newtheorem{lemma}{Lemma}
\newtheorem{notation}{Notation}

\newtheorem{proposition}{Proposition}
\newtheorem{remark}{Remark}

\numberwithin{equation}{section}
\input{tcilatex}

\begin{document}
\title[Good Rough Path Sequences]{Good Rough Path Sequences and Applications
to Anticipating \& Fractional Stochastic Calculus}
\author{Laure Coutin, Peter Friz, Nicolas Victoir}

\begin{abstract}
We consider anticipative Stratonovich stochastic differential equations
driven by some stochastic process (not necessarily a semi-martingale). No
adaptedness of initial point or vector fields is assumed. Under a simple
condition on the stochastic process, we show that the unique solution of the
above SDE understood in the rough path sense is actually a Stratonovich
solution. This condition is satisfied by the Brownian motion and the
fractional Brownian motion with Hurst parameter greater than $1/4$. As
application, we obtain rather flexible results such as support theorems,
large deviation principles and Wong-Zakai approximations for SDEs driven by
fractional Brownian Motion along anticipating vectorfields. In particular,
this unifies many results on anticipative SDEs.
\end{abstract}

\maketitle

\section{Introduction}

We fix a filtered probability space $\left( \Omega ,\mathcal{F},\mathbb{P}%
,\left( \mathcal{F}_{t}\right) _{0\leq t\leq 1}\right) $ satisfying the
usual conditions. It\^{o}'s theory tells us that there exists a unique
solution to the Stratonovich stochastic differential equation%
\begin{equation}
\left\{ 
\begin{array}{l}
dY_{t}=V_{0}(Y_{t})dt+\sum_{i=1}^{d}V_{i}(Y_{t})\circ dB_{t}^{i} \\ 
Y_{0}=y_{0},%
\end{array}%
\right.  \label{SDE}
\end{equation}%
where $B=\left( B_{t}^{1},\cdots ,B_{t}^{d}\right) _{0\leq t\leq 1}$ is a
standard $d$-dimensional $\left( \mathcal{F}_{t}\right) $-Brownian motion, $%
V_{0}$ a $C^{1}$ vector field on $\mathbb{R}^{d}$, $V_{1},\cdots ,V_{d}$
some $C^{2}$ vector fields on $\mathbb{R}^{d}$, and $y_{0}$ a $\mathcal{F}%
_{0}$-measurable random variable.

We remind the reader that, by definition (see \cite{Nu}), a process $z$ is
Stratonovich integrable with respect to a process $x$ if for all $t$ there
exists a random variable denoted $\int_{0}^{t}z_{u}\circ dx_{u}$ (the
Stratonovich integral of $z$ with respect to $x$ between time $0$ and $t$)
such that for all sequences $\left( D^{n}=\left( t_{i}^{n}\right)
_{i}\right) _{n\geq 0}$of subdivisions of $[0,t]$ such that $\left\vert
D^{n}\right\vert \rightarrow _{n\rightarrow \infty }0$, the following
convergence holds in probability%
\begin{equation*}
\sum_{i}\left( \frac{1}{t_{i+1}^{n}-t_{i}^{n}}%
\int_{t_{i}^{n}}^{t_{i+1}^{n}}z_{u}du\right) \left(
x_{t_{i+1}^{n}}-x_{t_{i}^{n}}\right) \rightarrow _{n\rightarrow \infty
}\int_{0}^{t}z_{u}\circ dx_{u}.
\end{equation*}%
Another way to express this is by introducing $x^{D}$ the $D$-linear
approximation of $x,$ where $D=\left( t_{i}\right) $ is a subdivision of $%
[0,1]$: 
\begin{equation}
x^{D}(t)=x_{t_{i}}+\frac{t-t_{i}}{t_{i+1}-t_{i}}\left(
x_{t_{i+1}}-x_{t_{i}}\right) \text{ if }t_{i}\leq t\leq t_{i+1}.
\label{Dlinear}
\end{equation}%
Then, $z$ is Stratonovich integrable with respect to $x$ \ if and only if
for all sequences of subdivisions $D^{n}$ which mesh size tends to $0$ and
for all $t\in \lbrack 0,1]$, $\int_{0}^{t}z_{u}dx^{D^{n}}(u)$ converges in
probability to $\int_{0}^{t}z_{u}\circ dx_{u}$.

Ocone and Pardoux \cite{OP} showed that there exists a unique solution to
equation (\ref{SDE}) even if the vector field $V_{0}$ and the initial
condition $y_{0}$ were allowed to be $\left( \mathcal{F}_{1}\right) $-random
variables. They did so relating the Skorokhod integral and the Stratonovich
one, and using Malliavin Calculus techniques. This solution has been studied
in various directions: existence and study of the density of $Y$ \cite%
{CFN,RS,SS}, a Freidlin-Wentzell's type theorem \cite{MNS}, results on the
support of the law of $\left\{ Y_{t},0\leq t\leq 1\right\} $ \cite{CFN2,MN},
approximation of $Y_{t}$ by some Euler's type schemes \cite{AK}... The case
where the vector fields $V_{1},\cdots ,V_{d}$ are allowed to be $\left( 
\mathcal{F}_{1}\right) $-measurable was dealt in \cite{KL,KLN}, under the
strong condition that the $\left( V_{i}\right) _{1\leq i\leq d}$ commute.
This is an application of the Doss-Sussman theorem. The latter says that if $%
V=\left( V_{1},\cdots ,V_{d}\right) $ are $d$ vector fields smooth enough
such that $\left[ V_{i},V_{j}\right] =0,$ for all $i,j\geq 1$, and if $V_{0}$
is another vector fields, then the map $\varphi _{\left( V_{0},V\right)
}^{y_{0}}$ which at a smooth path $x:\left[ 0,1\right] \rightarrow \mathbb{R}%
^{d}$ associates the path $y$ which is the solution of the differential
equation%
\begin{equation}
\left\{ 
\begin{array}{l}
dy_{t}=V_{0}(y_{t})dt+\sum_{i=1}^{d}V_{i}(y_{t})dx_{t}^{i} \\ 
y_{0}=y_{0},%
\end{array}%
\right.  \label{ODE}
\end{equation}%
is continuous when one equips the space of continuous functions with the
uniform topology. One can then define $\varphi _{\left( V_{0},V\right)
}^{y_{0}}$ on the whole space of continuous function, in particular $\varphi
_{\left( V_{0},V\right) }^{y_{0}}(B)$ is then well defined, and is
almost-surely the solution of the Stratonovich differential equation (\ref%
{SDE}). This remains true even if the vector fields and the initial
condition are allowed to be random.

Rough path theory can be seen as a major extension of the Doss-Sussman
result. One of the main thing to remember from this theory is that it is not 
$x$ which controls the differential equation (\ref{ODE}), but the lift of $x$
to a path in a Lie group lying over $\mathbb{R}^{d}$. The choice of the Lie
group depends on the roughness of $x$.

If $x$ is a $\mathbb{R}^{d}$-valued path of finite $p$-variation, $p\geq 1$,
one needs to lift $x$ to a path $\mathbf{x}$ with values in $G^{[p]}\left( 
\mathbb{R}^{d}\right) $, the free nilpotent Lie group of step $\left[ p%
\right] $ over $\mathbb{R}^{d}$. When $x$ is smooth, there exists a
canonical lift of $x$ to a path denoted $S(x)$ with values in $G^{[p]}\left( 
\mathbb{R}^{d}\right) $ ($S(x)$ is obtain from $x$ by computing the "first $%
[p]$"\ iterated integrals of $x$). If $x$ is a smooth path, then there
exists a solution $y$ to the differential equation (\ref{ODE}). We denote by 
$I_{\left( V_{0},V\right) }^{y_{0}}$ the map which at $S(x)$ associates $%
S(x\oplus y).$ Denoting $C\left( G^{[p]}\left( \mathbb{R}^{d}\right) \right) 
$ the set of continuous paths from $[0,1]$ into $G^{[p]}\left( \mathbb{R}%
^{d}\right) $, we see that $I_{\left( V_{0},V\right) }^{y_{0}}$ is a map
from a subset of $C\left( G^{[p]}\left( \mathbb{R}^{d}\right) \right) $ to $%
C\left( G^{[p]}\left( \mathbb{R}^{d}\oplus \mathbb{R}^{n}\right) \right) $.
Lyons showed that this map is (locally uniformly) continuous when one equips 
$C\left( G^{[p]}\left( \mathbb{R}^{d}\right) \right) $ and $C\left(
G^{[p]}\left( \mathbb{R}^{d}\oplus \mathbb{R}^{n}\right) \right) $ with a "$%
p $-variation distance". Hence one can define $I_{\left( V_{0},V\right)
}^{y_{0}}$ on the closure (in this $p$-variation topology) of the canonical
lift of smooth paths. This latter set is the set of geometric $p$-rough
paths.

\bigskip

\bigskip

\bigskip

In the case $x=B$, the Brownian motion being almost surely $1/p$-H\"{o}lder, 
$2<p<3$, one needs to lift $B$ to a process with values in $G^{2}(\mathbb{R}%
^{d})$ to obtain the solution in the rough path sense of equation (\ref{SDE}%
). This is equivalent to define its area process. The standard choice for
the area process of the Brownian motion is the L\'{e}vy area \cite{Levy,LQ},
although one could choose very different area processes \cite{LL}. Choosing
this area, we lift $B$ to a geometric $p$-rough path $\mathbf{B}$, $%
I_{\left( V_{0},V\right) }^{y_{0}}\left( \mathbf{B}\right) $ is then the
Stratonovich solution of the stochastic differential equation (\ref{SDE}),
together with its lift (i.e. here its area process).

Just as before, when the vector fields $V_{i}$ are almost surely "smooth
enough" and $y_{0}$ is almost surely finite, the It\^{o} map $I_{\left(
V_{0},V\right) }^{y_{0}}$ is still well defined and continuous almost
surely, and there is no problem at all of definition of $I_{\left(
V_{0},V\right) }^{y_{0}}\left( \mathbf{B}\right) $. Therefore, the theory of
rough path provides a meaning and a unique solution to the stochastic
differential equation (\ref{SDE}), even when the vector fields and the
initial condition depend on the whole Brownian path. Moreover, the
continuity of the It\^{o} map provides for free a Wong-Zakai theorem, and is
very well adapted to obtaining large deviation principles and support
theorems.

The only work not completely given for free by the theory of rough path is
to prove that the solution $\mathbf{y}$ of equation (\ref{SDE}) using the
rough path approach is actually solution of the Stratonovich differential
equation, i.e. that for all $t$,%
\begin{equation*}
y_{t}=y_{0}+\int_{0}^{t}V_{0}(y_{u})du+\sum_{i=1}^{d}%
\int_{0}^{t}V_{i}(y_{u})\circ dB_{u}^{i}.
\end{equation*}%
When the vector fields and the initial condition are deterministic, this is
usually proved using the standard Wong-Zakai theorem.

We provide in this general case here a solution typically in the spirit of
rough path, by separating neatly probability theory and differential
equation theory. We will show via a deterministic argument that to obtain
our result we only need to check that, if $D^{n}$ is a sequence of
subdivisions which steps tends to $0$, $\int_{0}^{.}B_{u}\otimes
dB_{u}^{D^{n}}$ and $\int_{0}^{.}B_{u}^{D^{n}}\otimes dB_{u}^{D^{n}}$
converges in an appropriate topology to $\int_{0}^{.}B_{u}\otimes \circ
dB_{u}$.

The paper is organized as follow:\ in the first section, we present quickly
the theory of rough path (for a more complete presentation, see \cite{Ly,LQ}
or \cite{Le}). The second section introduce the notion of good rough path
sequence and its properties. We will then show that $B^{n}$ defines a good
rough path sequence associated to $\mathbf{B}$, and this will imply that the
solution via rough path of equation (\ref{ODE}) with signal $\mathbf{B}$ is
indeed solution of the Stratonovich stochastic differential equation (\ref%
{SDE}).

We conclude with a few applications: an approximation/Wong-Zakai result, a
large deviation principle, and some remarks on the support theorem.

\section{Rough Paths}

By path we will always mean a continuous function from $[0,1]$ into a (Lie)\
group. If $x$ is such a path, $x_{s,t}$ is a notation for $x_{s}^{-1}.x_{t}$.

\subsection{Algebraic preliminaries}

We present the theory of rough paths, in the finite dimensional case purely
for simplicity. All arguments are valid in infinite dimension. We equip $%
\mathbb{R}^{d}$ with the Euclidean scalar product $\left\langle
.,.\right\rangle $ and the Euclidean norm $\left\vert x\right\vert
=\left\langle x,x\right\rangle ^{1/2}$. We denote by $\left( G^{n}(\mathbb{R}%
^{d}),\otimes \right) $ the free nilpotent of step $n$ over $\mathbb{R}^{d}$%
, which is imbedded in the tensor algebra $T^{n}(\mathbb{R}%
^{d})=\bigoplus_{k=0}^{n}\left( \mathbb{R}^{d}\right) ^{\otimes k}$. We
define a family of dilations on the group $G^{n}(\mathbb{R}^{d})$ by the
formula%
\begin{equation*}
\delta _{\lambda }(1,x_{1},\cdots ,x_{n})=(1,\lambda x_{1},\cdots ,\lambda
^{n}x_{n}),
\end{equation*}%
where $x_{i}\in \left( \mathbb{R}^{d}\right) ^{\otimes i},$ $(1,x_{1},\cdots
,x_{n})\in G^{n}(\mathbb{R}^{d})$ and $\lambda \in \mathbb{R}$. Inverse,
exponential and logarithm functions on $T^{n}(\mathbb{R}^{d})$ are defined
by mean of their power series \cite{Ly,Re}. We define on $\left( \mathbb{R}%
^{d}\right) ^{\otimes k}$ the Hilbert tensor scalar product and its norm%
\begin{eqnarray*}
\left\langle \sum_{i=1}^{n}x_{i}^{1}\otimes \cdots \otimes
x_{i}^{k},\sum_{j=1}^{n}y_{j}^{1}\otimes \cdots \otimes
y_{i}^{k}\right\rangle &=&\sum_{i,j}\left\langle
x_{i}^{1},y_{j}^{1}\right\rangle \cdots \left\langle
x_{i}^{k},y_{j}^{k}\right\rangle , \\
\left\vert x\right\vert &=&\left\langle x,x\right\rangle ^{1/2}.
\end{eqnarray*}%
This yields a family of comptible tensor norms on $\mathbb{R}^{d}$ and its
tensor product spaces. Since all finite dimensional norms are equivalent the
Hilbert structure of $\mathbb{R}^{d}$ was only used for convenience. In
fact, one can replace $\mathbb{R}^{d}$ by a Banach-space and deal with
suited tensor norms but this can be rather subtle, see \cite{LQ} and the
references therein.

\bigskip

For $(1,x_{1},\cdots ,x_{n})\in G^{n}(\mathbb{R}^{d})$, with $x_{i}\in
\left( \mathbb{R}^{d}\right) ^{\otimes i}$, we define%
\begin{equation*}
\left\Vert (1,x_{1},\cdots ,x_{n})\right\Vert =\max_{i=1,\cdots ,n}\left\{
\left( i!\left\vert x_{i}\right\vert \right) ^{1/i}\right\} .
\end{equation*}%
$\left\Vert .\right\Vert $ is then a symmetric sub-additive homogeneous norm
on $G^{n}(\mathbb{R}^{d})$ ($\left\Vert g\right\Vert =0$ iff $g=1$, and $%
\left\Vert \delta _{\lambda }g\right\Vert =\left\vert \lambda \right\vert
.\left\Vert g\right\Vert $ for all $\left( \lambda ,g\right) \in \mathbb{R}%
\times G^{n}(\mathbb{R}^{d})$, for all $g,h\in G^{n}(\mathbb{R}^{d})$, $%
\left\Vert g\otimes h\right\Vert \leq \left\Vert g\right\Vert +\left\Vert
h\right\Vert $, and $\left\Vert g^{-1}\right\Vert =\left\Vert g\right\Vert $%
). This homogeneous norm allows us to define on $G^{n}(\mathbb{R}^{d})$ a
left-invariant distance with the formula%
\begin{equation*}
d(g,h)=\left\Vert g^{-1}\otimes h\right\Vert \text{.}
\end{equation*}

From this distance we can define some distances on the space of continuous
paths from $[0,1]$ into $G^{[p]}(\mathbb{R}^{d})$ ($p$ is a fixed real
greater than or equal to $1$):

\begin{enumerate}
\item the $p$-variation distance%
\begin{equation*}
d_{p-var}(x,y)=d\left( x_{0},y_{0}\right) +\sup \left( \sum_{i}d\left(
x_{t_{i},t_{i+1}},y_{t_{i},t_{i+1}}\right) ^{p}\right) ^{1/p}
\end{equation*}%
where the supremum is over all subdivision $\left( t_{i}\right) _{i}$ of $%
[0,1]$. We also define $\left\Vert x\right\Vert _{p-var}=d_{p-var}\left(
1,x\right) $, and $\left\Vert x\right\Vert _{p-var}<\infty $ means that $x$
has finite $p$-variation.

\item modulus distances: We say that $\omega :\left\{ (s,t),0\leq s\leq
t\leq 1\right\} \rightarrow \mathbb{\ R}^{+}$ is a control if%
\begin{equation*}
\left\{ 
\begin{array}{l}
\omega \text{ is continuous.} \\ 
\omega \text{ is super-additive, i.e. }\forall s<t<u,\text{ }\omega
(s,t)+\omega (t,u)\leq \omega (t,u). \\ 
\omega \text{ is zero on the diagonal, i.e. }\omega (t,t)=0\text{ for all }%
t\in \lbrack 0,1]%
\end{array}%
\right.
\end{equation*}
\end{enumerate}

When $\omega $ is non-zero off the diagonal we introduce the distance $%
d_{\omega ,p}$ by%
\begin{equation*}
d_{\omega ,p}(x,y)=d\left( x_{0},y_{0}\right) +\sup_{0\leq s<t\leq 1}\frac{%
d\left( x_{s,t},y_{s,t}\right) }{\omega (s,t)^{1/p}}.
\end{equation*}%
The simplest example of such a control is given by $\omega (s,t)=t-s$. In
this case, $d_{\omega ,p}$ defines then a notion of $1/p$-H\"{o}lder
distance on $G^{n}(\mathbb{R}^{d})$-valued paths. In general, if $\left\Vert
x\right\Vert _{\omega ,p}=d_{\omega ,p}\left( 1,x\right) <\infty $, $\ $we
say that $x$ has finite $p$-variation controlled by $\omega $.

\subsection{$G^{n}(\mathbb{R}^{d})$-Valued Paths}

Let $x:[0,1]\rightarrow \mathbb{R}^{d}$ be a path of bounded variation and
define $S(x)$ to be the solution of the ordinary differential equation%
\begin{eqnarray*}
dS(x)_{t} &=&S(x)_{t}\otimes dx_{t} \\
S(x)_{0} &=&\exp (x_{0}),
\end{eqnarray*}%
where $\otimes $ is the multiplication of $T^{n}(\mathbb{R}^{d})$ and $\exp $
is the exponential on $T^{n}(\mathbb{R}^{d})$. $S(x)$ actually takes its
values in $G^{n}(\mathbb{R}^{d})$. Moreover, if the $1$-variation of $x$ is
controlled by $\omega $, then so is the $1$-variation of $S(x)$ \cite[%
Theorem 1]{Ly}.

\begin{example}
If $n=2$, $S(x)_{t}=\exp \left( x(t)+\frac{1}{2}\int_{0}^{t}x_{u}\otimes
dx_{u}-\frac{1}{2}\int_{0}^{t}dx_{u}\otimes x_{u}\right) $. The term $\frac{1%
}{2}\int_{0}^{t}x_{u}\otimes dx_{u}-\frac{1}{2}\int_{0}^{t}dx_{u}\otimes
x_{u}$ is the area between the line joining $x_{0}$ and $x_{t}$ and the path 
$\left( x_{u}\right) _{0\leq u\leq t}$.
\end{example}

\begin{definition}
A path $\mathbf{x}:[0,1]\rightarrow G^{[p]}(\mathbb{R}^{d})$ is a geometric $%
p$-rough path if there exists a sequence of paths of bounded variation $%
\left( x_{n}\right) _{n\in \mathbb{N}}$ such that%
\begin{equation*}
\lim_{n\rightarrow \infty }d_{p-var}\left( S\left( x_{n}\right) ,\mathbf{x}%
\right) =0.
\end{equation*}
\end{definition}

\begin{notation}
Let $C^{0,p-var}\left( G^{[p]}(\mathbb{R}^{d})\right) $ denote the space of
geometric $p$-rough paths, and $C^{0,\omega ,p}\left( G^{[p]}(\mathbb{R}%
^{d})\right) $ the set of paths $\mathbf{x}:[0,1]\rightarrow G^{[p]}(\mathbb{%
R}^{d})$ for which there exists a sequence of paths $\left( x_{n}\right)
_{n\in \mathbb{N}}$ of bounded variation such that%
\begin{equation*}
\lim_{n\rightarrow \infty }d_{\omega ,p}\left( S\left( x_{n}\right) ,\mathbf{%
x}\right) =0.
\end{equation*}
\end{notation}

In particular, any $G^{[p]}(\mathbb{R}^{d})$-valued path with finite $q$%
-variation, $q<p$, is a geometric $p$-rough path, and a geometric $p$-rough
path has finite $p$-variation. We refer to \cite{FV} for a "Polish" study of
the space of geometric rough paths.

If $\mathbf{x}$ is a geometric $p$-rough path, we denote by $\mathbf{x}%
_{s,t}^{i}$ the projection of $\mathbf{x}$ onto $\left( \mathbb{R}%
^{d}\right) ^{\otimes i}$. We also say that $\mathbf{x}$ lies above $\mathbf{%
x}^{1}$.

\subsection{The It\^{o} map}

We fix $p\geq 1$, and a control $\omega $.

When $A$ is vector field on $\mathbb{R}^{d}$, we denote by $d^{k}A$ its $%
k^{th}$ derivative (with the convention $d^{0}A=A$), and its $\gamma $%
-Lipschitz norm by%
\begin{equation*}
\left\Vert A\right\Vert _{Lip\left( \gamma \right) }=\max \left\{
\max_{k=0,\cdots ,[\gamma ]}\left\Vert d^{k}A\right\Vert _{\infty
},\left\Vert d^{[\gamma ]}A\right\Vert _{\gamma -[\gamma ]}\right\} ,
\end{equation*}%
where $\left\Vert .\right\Vert _{\infty }$ is the sup norm and $\left\Vert
.\right\Vert _{\beta }$ the $\beta $-H\"{o}lder norm, $0\leq \beta <1$. If $%
\left\Vert A\right\Vert _{Lip\left( \gamma \right) }<\infty ,$ we say that $%
A $ is a $Lip\left( \gamma \right) $-vector fields on $\mathbb{R}^{d}$.

We consider $V=\left( V_{1},\cdots ,V_{d}\right) $, where the $V_{i}$ are $%
Lip\left( p+\varepsilon \right) $-vector fields on $\mathbb{R}^{n}$, $%
\varepsilon >0$. $V$ can be identified with a linear map from $\mathbb{R}%
^{d} $ into $Lip\left( p+\varepsilon \right) $-vector fields on $\mathbb{R}%
^{n},$ 
\begin{equation*}
V(y)(dx^{1},\cdots ,dx^{d})=\sum_{i=1}^{d}V_{i}(y)dx^{i}.
\end{equation*}%
For a $\mathbb{R}^{d}$-valued path $x$ of bounded variation, we define $y$
to be the solution of the ordinary differential equation%
\begin{equation}
\left\{ 
\begin{array}{l}
dy_{t}=V(y_{t})dx_{t} \\ 
y_{0}=y_{0}.%
\end{array}%
\right.  \label{ODE2}
\end{equation}%
Lifting $x$ and $\left( x\oplus y\right) $ to $S(x)$ and $S(x\oplus y)$
(their canonical lift to paths with values in the free nilpotent group of
step $p$), we consider the map which at $S(x)$ associates $S(x\oplus y)$. We
denote it $I_{y_{0},V}$. We refer to \cite{Ly,LQ} for the following theorem:

\begin{theorem}[Universal Limit Theorem (Lyons)]
\label{ULT}The map $I_{y_{0},V}$ is continuous from $\left( C^{0,\omega
,p}\left( G^{[p]}(\mathbb{R}^{d})\right) ,d_{\omega ,p}\right) $ into $%
\left( C^{0,\omega ,p}\left( G^{[p]}(\mathbb{R}^{d}\oplus \mathbb{R}%
^{n})\right) ,d_{\omega ,p}\right) $.
\end{theorem}

Let $x_{n}$ be a sequence of path of bounded variation such that $S\left(
x_{n}\right) $ converges in the $d_{\omega ,p}$-topology to a geometric $p$%
-rough path $\mathbf{x}$, and define $y_{n}$ to be the solution of the
differential equation (\ref{ODE2}), where $x$ is replaced by $x_{n}$. Then,
the Universal Limit Theorem says that $S(x_{n}\oplus y_{n})$ converges in
the $d_{\omega ,p}$-topology to a geometric $p$-rough path $\mathbf{z}$. We
say that $\mathbf{y}$\textbf{,} the projection of $\mathbf{z}$ onto $G^{[p]}(%
\mathbb{R}^{n})$ is the solution of the rough differential equation%
\begin{equation*}
d\mathbf{y}=V(\mathbf{y})d\mathbf{x}
\end{equation*}%
with initial condition $y_{0}$. It is interesting to observe that Lyons'
estimates actually give that for all $R>0$ and sequence $\mathbf{x}_{n}\in
C^{0,\omega ,p}\left( G^{[p]}(\mathbb{R}^{d}\oplus \mathbb{R}^{n})\right) $
converging to $\mathbf{x}\in C^{0,\omega ,p}\left( G^{[p]}(\mathbb{R}%
^{d}\oplus \mathbb{R}^{n})\right) $ in the $d_{\omega ,p}$-topology,%
\begin{equation}
\sup_{\substack{ \left\vert y_{0}\right\vert \leq R  \\ \left\Vert
V\right\Vert _{Lip\left( p+\varepsilon \right) }\leq R}}d_{\omega ,p}\left(
I_{y_{0},V}\left( \mathbf{x}_{n}\right) ,I_{y_{0},V}\left( \mathbf{x}\right)
\right) \rightarrow _{n\rightarrow \infty }0.
\label{continuity_uniform_over_bounded_set}
\end{equation}%
To see this, one first observes that the continuity of integration along of
a one-form is uniform over the set of one-forms with Lipschitz norm bounded
by a given $R$. Then, the path $I_{y_{0},V}\left( \mathbf{x}\right) $ over
small time is the fixed point of a map (integrating along a one-form) which
is a contraction. Reading the estimate in \cite{LQ}, one sees that this map
is uniformly, over the set $\left\{ \left\vert y_{0}\right\vert \leq
R;\left\Vert V\right\Vert _{Lip\left( p+\varepsilon \right) }\leq R\right\}
, $ a contraction with parameter strictly less than $1$.

The next theorem was also proved in \cite{LQ}, and deals with the continuity
of the flow.

\begin{theorem}
\label{ITOIC}If $\left( y_{0}^{n}\right) _{n}$ is a $\mathbb{R}^{d}$-valued
sequence converging to $y_{0}$, then for all $R>0,$%
\begin{equation*}
\sup_{_{\substack{ \left\Vert \mathbf{x}\right\Vert _{\omega ,p}\leq R  \\ %
\left\Vert V\right\Vert _{Lip\left( p+\varepsilon \right) }\leq R}}%
}d_{\omega ,p}\left( I_{y_{0}^{n},V}\left( \mathbf{x}\right)
,I_{y_{0},V}\left( \mathbf{x}\right) \right) \rightarrow _{n\rightarrow
\infty }0.
\end{equation*}
\end{theorem}

The next theorem shows that the It\^{o} map it continuous when one varies
the vector fields defining the differential equation. It does not seem to
have appeared anywhere, despite that its proof does not involve any new
ideas.

\begin{theorem}
\label{ITOVF}Let $V^{n}=\left( V_{1}^{n},\cdots ,V_{d}^{n}\right) $ be a
sequence of $d$ $Lip\left( p+\varepsilon \right) $-vector fields on $\mathbb{%
R}^{n}$ such that%
\begin{equation*}
\lim_{n\rightarrow \infty }\max_{i}\left\Vert V_{i}^{n}-V_{i}\right\Vert
_{Lip\left( p+\varepsilon \right) }=0.
\end{equation*}%
Then, for all $R>0,$%
\begin{equation*}
\lim_{n\rightarrow \infty }\sup_{\substack{ \left\Vert \mathbf{x}\right\Vert
_{\omega ,p}\leq R  \\ \left\vert y_{o}\right\vert \leq R}}d_{\omega
,p}\left( I_{y_{0},V^{n}}\left( \mathbf{x}\right) ,I_{y_{0},V}\left( \mathbf{%
x}\right) \right) =0.
\end{equation*}
\end{theorem}

\begin{proof}
We use the notations of \cite{Ly}. First consider the $Lip\left(
p+\varepsilon -1\right) $-one-forms $\theta _{i}:\mathbb{R}^{d}\rightarrow
Hom\left( \mathbb{R}^{d},\mathbb{R}^{n}\right) $, $i\in \mathbb{N\cup }%
\left\{ \infty \right\} $. We assume that $\theta _{n}$ converges to $\theta
_{\infty }$ in the $\left( p+\varepsilon -1\right) $-Lipschitz topology when 
$n\rightarrow \infty $. For $n\in \mathbb{N\cup }\left\{ \infty \right\} $, $%
\int \theta _{n}\left( \mathbf{x}\right) d\mathbf{x}$ is the unique rough
path associated to the almost multiplicative functional%
\begin{equation*}
\left( Z\left( \mathbf{x}\right) _{s,t}^{n}\right) ^{i}=\sum_{l_{1},\cdots
,l_{i}=0}^{[p]-1}d^{l_{1}}\theta _{n}\left( \mathbf{x}_{s}^{1}\right)
\otimes \cdots \otimes d^{l_{i}}\theta _{n}\left( \mathbf{x}_{s}^{1}\right)
\left( \sum_{\pi \in \Pi _{\left( l_{1},\cdots ,l_{i}\right) }}\pi \left( 
\mathbf{x}_{s,t}^{i+l_{1}+\cdots +l_{i}}\right) \right)
\end{equation*}%
(see \cite{Ly} for the definition of $\Pi $). It is obvious that%
\begin{equation*}
\sup_{\left\Vert \mathbf{x}\right\Vert _{\omega ,p}\leq R}\max_{i}\frac{%
\left\vert \left( Z\left( \mathbf{x}\right) _{s,t}^{n}\right) ^{i}-\left(
Z\left( \mathbf{x}\right) _{s,t}^{\infty }\right) ^{i}\right\vert }{\omega
\left( s,t\right) ^{i/p}}\rightarrow _{n\rightarrow \infty }0,
\end{equation*}%
which implies by theorem 3.1.2 in \cite{Ly} that%
\begin{equation*}
\sup_{_{\left\Vert \mathbf{x}\right\Vert _{\omega ,p}\leq R}}d_{\omega
,p}\left( \int \theta _{n}\left( \mathbf{x}\right) d\mathbf{x,}\int \theta
_{\infty }\left( \mathbf{x}\right) d\mathbf{x}\right) \rightarrow
_{n\rightarrow \infty }0.
\end{equation*}%
Now consider the Picard iteration sequence $\left( \mathbf{z}_{m}^{n}\right)
_{m\geq 0}$ introduced in \cite[formula (4.10)]{Ly} to construct to $%
I_{y_{0},V^{n}}\left( \mathbf{x}\right) $:%
\begin{eqnarray*}
\mathbf{z}_{0}^{n} &=&\left( 0,y_{0}\right) \\
\mathbf{z}_{m+1}^{n} &=&\int h_{n}\left( \mathbf{z}_{m}^{n}\right) d\mathbf{z%
}_{m}^{n},
\end{eqnarray*}%
where $h_{n}$ is the one-form defined by the formula $h_{n}\left( x,y\right)
\left( dX,dY\right) =\left( dX,V\left( y\right) dX\right) $ (by $V(y)dX$ we
mean $\sum V_{i}(y)dX^{i}$). Lyons proved that%
\begin{equation*}
\lim_{m\rightarrow \infty }\sup_{\substack{ \left\Vert \mathbf{x}\right\Vert
_{\omega ,p}\leq R  \\ \left\vert y_{o}\right\vert \leq R}}\sup_{n}d_{\omega
,p}\left( \mathbf{z}_{m}^{n},\mathbf{z}_{\infty }^{n}\right) =0;
\end{equation*}%
we have the supremum over all $n$ here because the $p+\varepsilon $%
-Lipschitz norm of the $V^{n}$ are uniformly bounded in $n$. Moreover, we
have just seen that for all fixed $m,$ $\lim_{n\rightarrow \infty
}\sup_{\left\vert y_{0}\right\vert \leq R}d_{\omega ,p}\left( \mathbf{z}%
_{m}^{n},\mathbf{z}_{m}^{\infty }\right) =0$. Therefore, with a $%
3\varepsilon $-type argument, we obtain our theorem.
\end{proof}

The three previous theorems actually gives that the map%
\begin{equation*}
\left( y_{0},V,\mathbf{x}\right) \rightarrow I_{y_{0},V}\left( \mathbf{x}%
\right)
\end{equation*}%
is continuous in the product topology $\mathbb{R}^{n}\times \left(
Lip(p+\epsilon )\text{ on }\mathbb{R}^{n}\right) ^{d}\times C^{0,\omega
,p}\left( G^{[p]}(\mathbb{R}^{d})\right) .$

In the reminder of this section, $p$ is a real in $[2,3)$. Now consider a $%
Lip\left( 1+\varepsilon \right) $-vector field on $\mathbb{R}^{n}$ denoted $%
V_{0}$, and for a path $x$ of bounded variation, we consider $y$ to be the
solution of%
\begin{equation}
\left\{ 
\begin{array}{l}
dy_{t}=V_{0}(y_{t})dt+V(y_{t})dx_{t} \\ 
y_{t}=y_{0},%
\end{array}%
\right.  \label{ODE1}
\end{equation}%
Lifting $x$ and $\left( x\oplus y\right) $ to $S(x)$ and $S\left( x\oplus
y\right) $ we consider the map which at $S(x)$ associates $S\left(
x,y\right) $. We denote it $I_{y_{0},\left( V_{0},V\right) }$. The following
extension of the Universal Limit Theorem was obtained in \cite{LV}.

\begin{theorem}
\label{ULT copy(1)}The map $I_{y_{0},\left( V_{0},V\right) }$ is continuous
from $\left( C^{0,\omega ,p}\left( G^{[p]}(\mathbb{R}^{d})\right) ,d_{\omega
,p}\right) $ into $\left( C^{0,\omega ,p}\left( G^{[p]}(\mathbb{R}^{d}\oplus 
\mathbb{R}^{n})\right) ,d_{\omega ,p}\right) $. More precisely, for all $R>0$%
, 
\begin{equation}
\sup_{\substack{ \left\vert y_{0}\right\vert \leq R  \\ \left\Vert
V\right\Vert _{Lip\left( p+\varepsilon \right) }\leq R}}d_{\omega ,p}\left(
I_{y_{0},\left( V_{0},V\right) }\left( \mathbf{x}_{n}\right)
,I_{y_{0},\left( V_{0},V\right) }\left( \mathbf{x}\right) \right)
\rightarrow _{n\rightarrow \infty }0.
\label{continuity_uniform_over_bounded_set2}
\end{equation}
\end{theorem}

Let $x_{n}$ be a sequence of paths of bounded variation such that $S\left(
x_{n}\right) $ converges in the $d_{\omega ,p}$-topology to a geometric $p$%
-rough path $\mathbf{x}$, and define $y_{n}$ to be the solution of equation (%
\ref{ODE1}) replacing $x$ by $x_{n}$. Then, the previous theorem says that $%
S(x_{n}\oplus y_{n})$ converges in the $d_{\omega ,p}$-topology to a
geometric $p$-rough path $\mathbf{z}$. We say that $\mathbf{y}$, the
projection of $\mathbf{z}$ onto $G^{[p]}\left( \mathbb{R}^{n}\right) $ is
the solution of the rough differential equation%
\begin{equation*}
d\mathbf{y}_{t}=V_{0}(\mathbf{y}_{t})dt+V(\mathbf{y}_{t})d\mathbf{x}_{t}
\end{equation*}%
with initial condition $y_{0}$.

We obtain, as before, the following two theorems:

\begin{theorem}
\label{ITOIC copy(1)}If $\left( y_{0}^{n}\right) _{n}$ is a $\mathbb{R}^{d}$%
-valued sequence converging to $y_{0}$, then 
\begin{equation*}
\sup_{\substack{ _{\substack{ \left\Vert \mathbf{x}\right\Vert _{\omega
,p}\leq R  \\ \left\Vert V\right\Vert _{Lip\left( p+\varepsilon \right)
}\leq R}}  \\ \left\Vert V_{0}\right\Vert _{Lip\left( 1+\varepsilon \right)
}\leq R }}d_{\omega ,p}\left( I_{y_{0}^{n},\left( V_{0},V\right) }\left( 
\mathbf{x}\right) ,I_{y_{0},\left( V_{0},V\right) }\left( \mathbf{x}\right)
\right) \rightarrow _{n\rightarrow \infty }0.
\end{equation*}
\end{theorem}

\begin{theorem}
\label{ITOVF copy(1)}Let $\left( V^{n}=\left( V_{1}^{n},\cdots
,V_{d}^{n}\right) \right) _{n\geq 0}$ be a sequence of $d$ $Lip\left(
p+\varepsilon \right) $-vector fields on $\mathbb{R}^{n}$ and $\left(
V_{0}^{n}\right) _{n\geq 0}$ a sequence of $Lip\left( 1+\varepsilon \right) $%
-vector fields on $\mathbb{R}^{n}$, such that%
\begin{equation*}
\lim_{n\rightarrow \infty }\max \left\{ \left\Vert
V_{0}^{n}-V_{0}\right\Vert _{Lip\left( 1+\varepsilon \right) },\max_{1\leq
i\leq d}\left\Vert V_{i}^{n}-V_{i}\right\Vert _{Lip\left( p+\varepsilon
\right) }\right\} =0.
\end{equation*}%
Then, if $\mathbf{x\in }C^{0,\omega ,p}\left( G^{[p]}(\mathbb{R}^{d})\right) 
$,%
\begin{equation*}
\lim_{n\rightarrow \infty }\sup_{\substack{ \left\Vert \mathbf{x}\right\Vert
_{\omega ,p}\leq R  \\ \left\vert y_{o}\right\vert \leq R}}d_{\omega
,p}\left( I_{y_{0},\left( V_{0}^{n},V^{n}\right) }\left( \mathbf{x}\right)
,I_{y_{0},\left( V_{0},V\right) }\left( \mathbf{x}\right) \right) =0.
\end{equation*}
\end{theorem}

\subsection{Solving Anticipative Stochastic Differential Equations Via Rough
Paths\label{ANT_VIA_ROUGH}}

We fix a $p\in (2,3)$ and, for simplicity, the control $\omega (s,t)=t-s$
(i.e. we deal with H\"{o}lder topologies), although we could have been more
general and have considered a wide class of controls as in \cite{FV1} (i.e.
we could have consider modulus type topologies). We define $\mathbf{B}$ the
Stratonovich lift to a geometric $p$-rough path of the Brownian motion $B$
with the formula%
\begin{equation*}
\mathbf{B}_{t}=\left( 1,B_{t},\int_{0}^{t}B_{u}\otimes \circ dB_{u}\right) .
\end{equation*}%
$\mathbf{B}$ is a $G^{2}\left( \mathbb{R}^{d}\right) $-valued path, and
almost surely, $\left\Vert \mathbf{B}\right\Vert _{\omega ,p}<\infty $.

Consider $V_{0}$ a random vector field on $\mathbb{R}^{n}$ almost surely in $%
Lip\left( 1+\varepsilon \right) ,$ i.e. a measurable map from%
\begin{equation*}
V_{0}:\Omega \times \mathbb{R}^{n}\rightarrow \mathbb{R}^{n}
\end{equation*}%
such that $V_{0}(\omega ,.)\in Lip\left( 1+\varepsilon \right) $ for $\omega 
$ in a set a full measure, and $V=\left( V_{1},\cdots ,V_{d}\right) ,$ where 
$V_{1},\cdots ,V_{d}$ are random vector fields on $\mathbb{R}^{n}$ almost
surely in $Lip\left( 2+\varepsilon \right) $, and a random variable $%
y_{0}\in \mathbb{R}^{n}$ finite almost surely.

$I_{y_{0},\left( V_{0},V\right) }\left( \mathbf{B}\right) $ is then almost
surely well defined, and its projection onto $G^{2}\left( \mathbb{R}%
^{n}\right) $ is the solution, in the rough path sense, of the anticipative
stochastic differential equation%
\begin{equation}
d\mathbf{y}_{t}=V_{0}(\mathbf{y}_{t})dt+V(\mathbf{y}_{t})d\mathbf{B}_{t}
\label{RDE}
\end{equation}%
with initial condition $y_{0}$.

The next section introduces the notion of good rough paths sequence, and its
properties. Showing that linear approximations of Brownian motion form good
rough path sequences (in some sense that will be precise later on) will
prove that $\mathbf{y}^{1}$ is solution of the anticipative Stratonovich
stochastic differential equation (\ref{SDE}). In particular, the solution
that we construct coincides with the one constructed in \cite{OP}.

\section{Good Rough Path Sequence}

\subsection{Definitions}

We fix a parameter $p>2$, and a control $\omega $. $\mathbb{R}^{d}$ and $%
\widetilde{\mathbb{R}^{d}}$ will denote two identical copies of $\mathbb{R}%
^{d}$.

Let $p>2$, and $q$ such that $1/p+1/q>1$. We consider $x$ and $y$ two $%
\mathbb{R}^{d}$-valued paths of bounded variation. We let $\mathbf{y=}S(y)$
to be the canonical lift of $y$ to a $G^{[p]}\left( \mathbb{R}^{d}\right) $%
-valued path. We let 
\begin{equation}
S^{\prime }(x,S(y)):=S(x\oplus y)  \label{SprimeLemma}
\end{equation}%
be the canonical lift of $x\oplus y$ to a $G^{[p]}\left( \mathbb{R}^{d}%
\mathbb{\oplus R}^{d}\right) $-valued path and 
\begin{equation}
S^{\prime \prime }(S(x)):=S^{\prime }(x,S(x))=S(x\oplus x)
\label{Sprimeprime}
\end{equation}
be the canonical lift of $y\oplus y$ to A $G^{[p]}\left( \mathbb{R}^{d}%
\mathbb{\oplus R}^{d}\right) $-valued path.

\begin{proposition}
$\label{contS'}$Let $x$ be a $R^{d}$-valued path of finite $q$-variation,
and $\mathbf{y}$ a $G^{[p]}\left( \mathbb{R}^{d}\right) $-valued path of
finite $p$-variation. Let $\left( x_{n},y_{n}\right) $ be a sequence of $%
\mathbb{R}^{d}\mathbb{\oplus R}^{d}$-valued path such that $d_{\omega
,p}\left( x_{n},x\right) \rightarrow _{n\rightarrow \infty }0$ and $%
d_{\omega ,p}\left( S(y_{n}),\mathbf{y}\right) \rightarrow _{n\rightarrow
\infty }0$. Then, \newline
(i) $S^{\prime }(x_{n},S(y_{n}))$ converges in the $d_{\omega ,p}$-topology,
and the limit is independent of the choice of the sequence $\left(
x_{n},y_{n}\right) $. We denote this limit element $S^{\prime }\left( x,%
\mathbf{y}\right) .$\newline
(ii) $S^{^{\prime \prime }}\left( S(y_{n}),S(y_{n})\right) $ converges in $%
d_{\omega ,p}$-topology, and the limit is independent of the choice of the
sequence $y_{n}$. We denote this elememt $S^{\prime \prime }\left( \mathbf{y}%
,\mathbf{y}\right) $.
\end{proposition}

\begin{proof}
This is simply obtained using theorem 3.1.2 in \cite{Ly}, which says that
the procedure which at an almost multiplicative functional associates a
rough path is continuous, and we leave the details to the reader.
\end{proof}

\begin{example}
If $2\leq p<3$, 
\begin{eqnarray*}
&&S^{\prime }(x,\mathbf{y})_{t} \\
&=&\left( 1,x_{t}\oplus \mathbf{y}_{t}^{1},\int_{0}^{t}x_{u}\otimes
dx_{u}\oplus \int_{0}^{t}x_{u}\otimes d\mathbf{y}_{u}^{1}\oplus \int_{0}^{t}%
\mathbf{y}_{u}^{1}\otimes dx_{u}\oplus \mathbf{y}_{t}^{2}\right) .
\end{eqnarray*}%
The three integrals are well defined Young's integrals \cite{Yo}.
\end{example}

We introduce the notion of a good $p$-rough path sequence.

\begin{definition}
Let $\left( x_{n}\right) _{n}$ be a sequence of $\mathbb{R}^{d}$-valued
paths of bounded variation, and $\mathbf{x}$ a geometric $p$-rough path. We
say that $\left( x_{n}\right) _{n\in \mathbb{N}}$ is a good $p$-rough path
sequence (associated to $\mathbf{x}$) (for the control $\omega $) if and
only if 
\begin{equation}
\lim_{n\rightarrow \infty }d_{\omega ,p}(S^{\prime }(x_{n},\mathbf{x}%
),S^{\prime \prime }(\mathbf{x}))=0.  \label{GOOD}
\end{equation}
\end{definition}

In particular, if $\left( x_{n}\right) _{n}$ is a good $p$-rough path
sequence associated to $\mathbf{x}$, for the control $\omega $, $x_{n}$
converges to $\mathbf{x}$ in the topology induced by $d_{\omega ,p}$.

\begin{proposition}
\label{good2}Assume $2\leq p<3$. The sequence $\left( x(n)\right) _{n}$ of
paths of bounded variation is good rough path sequence associated to $%
\mathbf{x,}$ for the control $\omega $, if and only if \newline
\begin{eqnarray*}
\lim_{n\rightarrow \infty }\sup_{0\leq s<t\leq 1}\frac{\left\vert x(n)_{s,t}-%
\mathbf{x}_{s,t}^{1}\right\vert }{\omega (s,t)^{1/p}} &=&0, \\
\lim_{n\rightarrow \infty }\sup_{0\leq s<t\leq 1}\frac{\left\vert
\int_{s}^{t}x(n)_{s,u}\otimes dx(n)_{u}-\mathbf{x}_{s,t}^{2}\right\vert }{%
\omega (s,t)^{2/p}} &=&0, \\
\lim_{n\rightarrow \infty }\sup_{0\leq s<t\leq 1}\frac{\left\vert
\int_{s}^{t}\mathbf{x}_{s,u}^{1}\otimes dx(n)_{u}-\mathbf{x}%
_{s,t}^{2}\right\vert }{\omega (s,t)^{2/p}} &=&0.
\end{eqnarray*}
\end{proposition}

\begin{proof}
$d_{\omega ,p}\left( S^{\prime }\left( x(n),\mathbf{x}\right) ,S^{\prime
\prime }\left( \mathbf{x},\mathbf{x}\right) \right) $ is bounded by a
constant times%
\begin{equation*}
\max \left\{ A_{1}^{n},\sqrt{A_{2}^{n}},\sqrt{A_{3}^{n}},\sqrt{A_{4}^{n}}%
\right\} ,
\end{equation*}
where%
\begin{eqnarray*}
A_{1}^{n} &=&\sup_{0\leq s<t\leq 1}\frac{\left\vert \mathbf{x}%
_{s,t}^{1}-x(n)_{s,t}\right\vert }{\omega (s,t)^{1/p}}, \\
A_{2}^{n} &=&\sup_{0\leq s<t\leq 1}\frac{\left\vert
\int_{s}^{t}x(n)_{s,u}\otimes dx(n)_{u}-\mathbf{x}_{s,t}^{2}\right\vert }{%
\omega (s,t)^{2/p}}, \\
A_{3}^{n} &=&\sup_{0\leq s<t\leq 1}\frac{\left\vert
\int_{s}^{t}x(n)_{s,u}\otimes d\mathbf{x}_{u}^{1}-\mathbf{x}%
_{s,t}^{2}\right\vert }{\omega (s,t)^{2/p}}, \\
A_{4}^{n} &=&\sup_{0\leq s<t\leq 1}\frac{\left\vert \int_{s}^{t}\mathbf{x}%
_{s,u}^{1}\otimes dx(n)_{u}-\mathbf{x}_{s,t}^{2}\right\vert }{\omega
(s,t)^{2/p}}.
\end{eqnarray*}%
Let $\pi $ be the linear operator from $\mathbb{R}^{d}\otimes \mathbb{R}^{d}$
onto itself defined by $\pi \left( x\otimes y\right) =\pi \left( y\otimes
x\right) $. Observe that $\pi \left( \mathbf{x}_{s,t}^{2}\right) =\left( 
\mathbf{x}_{s,t}^{1}\right) ^{\otimes 2}-\mathbf{x}_{s,t}^{2}$. Also, we
have for all $z\in \mathbb{R}^{d}\otimes \mathbb{R}^{d}$, $\left\vert
z\right\vert =\left\vert \pi \left( z\right) \right\vert .$ Using this
property and an integration by part, we see that%
\begin{eqnarray*}
&&\left\vert \int_{s}^{t}\mathbf{x}_{s,u}^{1}\otimes dx(n)_{u}-\mathbf{x}%
_{s,t}^{2}\right\vert \\
&=&\left\vert \mathbf{x}_{s,t}^{1}\otimes x(n)_{s,t}-\int_{s}^{t}d\mathbf{x}%
_{s,u}^{1}\otimes x(n)_{u}-\mathbf{x}_{s,t}^{2}\right\vert \\
&=&\left\vert \mathbf{x}_{s,t}^{1}\otimes \left( x(n)_{s,t}-\mathbf{x}%
_{s,t}^{1}\right) -\int_{s}^{t}d\mathbf{x}_{s,u}^{1}\otimes x(n)_{u}+\left(
\left( \mathbf{x}_{s,t}^{1}\right) ^{\otimes 2}-\mathbf{x}_{s,t}^{2}\right)
\right\vert \\
&\leq &\left\vert \mathbf{x}_{s,t}^{1}\right\vert .\left\vert x(n)_{s,t}-%
\mathbf{x}_{s,t}^{1}\right\vert +\left\vert \mathbf{x}_{s,t}^{2}-%
\int_{s}^{t}x(n)_{s,u}\otimes d\mathbf{x}_{u}^{1}\right\vert .
\end{eqnarray*}%
Hence, $A_{4}^{n}\leq \left\Vert \mathbf{x}\right\Vert _{\omega
,p}.A_{1}^{n}+A_{3}^{n}$, which proves the proposition.
\end{proof}

By definition of a geometric $p$-rough path, there always exists a sequence
of smooth $x_{n}$ such that $S(x_{n})$ converges to $\mathbf{x}$ in the
topology induced by $d_{\omega ,p}$. However, this does not imply that
equation (\ref{GOOD}) holds (or equivalently for $[p]=2:$ that the
conditions in proposition \ref{good2} hold). We now give an example of a
geometric rough path $\mathbf{x}$ for which there exists no good sequence
associated to it.

\begin{example}
Consider%
\begin{eqnarray*}
\mathbf{x}_{t} &=&\exp (t[e_{1},e_{2}]) \\
&=&\exp (t(e_{1}\otimes e_{2}-e_{2}\otimes e_{1}))
\end{eqnarray*}%
where $e_{1},e_{2}$ is a basis of $\mathbb{R}^{2}$. Note that $\mathbf{x}$
is a geometric $p$-rough path and%
\begin{equation*}
\mathbf{x}_{t}^{1}\equiv 0\in \mathbb{R}^{2},\text{ \ \ \ }x_{t}^{2}\equiv
t[e_{1},e_{2}]\in \mathbb{R}^{2}\otimes \mathbb{R}^{2}.
\end{equation*}
If $x_{n}$ is a sequence of smooth paths in $\mathbb{R}^{2}$, then $%
\int_{0}^{t}\mathbf{x}^{1}\otimes dx_{n}\equiv 0$ trivially converges to $%
0\in \mathbb{R}^{2}\otimes \mathbb{R}^{2}.$ On the other hand, if $x_{n}$
was a good $p$-rough path sequence associated to $\mathbf{x}$ it should
converge to $t[e_{1,}e_{2}]$
\end{example}

\subsection{Stability of Good Rough Path Sequences}

Nonetheless, there are some good news. Good rough path sequences are stable
under integration, as shows the following theorem. We fix a control $\omega $%
.

\begin{theorem}
\label{staygood}Let $\left( x_{n}\right) _{n\in \mathbb{N}}$ denote a good $%
p $-rough path sequence associated to $\mathbf{x}$ for the control $\omega $%
, and $\theta :\mathbb{R}^{d}\rightarrow Hom(\mathbb{R}^{d},\mathbb{R}^{n})$
be a $Lip\left( p-1+\varepsilon \right) $ one-form, with $\varepsilon >0$.
Then $\left( \int \theta (x_{n})dx_{n}\right) _{n}$ is a good rough path
sequence associated to $\int \theta (\mathbf{x})d\mathbf{x}$ for the control 
$\omega $. Moreover, 
\begin{equation}
\lim_{n\rightarrow \infty }S\left( \int \theta (\mathbf{x}^{1})dx_{n}\right)
=\int \theta (\mathbf{x})d\mathbf{x}.  \label{givesStrat}
\end{equation}
\end{theorem}

\begin{proof}
Consider the $Lip\left( p-1+\varepsilon \right) $ one-form%
\begin{eqnarray*}
\widetilde{\theta } &:&\mathbb{R}^{d}\oplus \mathbb{R}^{d}\rightarrow Hom(%
\mathbb{R}^{d}\oplus \mathbb{R}^{d},\mathbb{R}^{n}\oplus \mathbb{R}^{n}) \\
((x,y),(dX,dY)) &\rightarrow &\theta (x)dX\oplus \theta (y)dY.
\end{eqnarray*}%
All the limits considered in this proof are to be understood to be in the
topology induced by $d_{\omega ,p}$. By the continuity of the integral \cite%
{Ly,LQ} and of the operator $S^{\prime }$ (by proposition \ref{contS'}), 
\begin{eqnarray*}
S^{\prime }\left( \int \theta (x_{n})dx_{n},\int \theta (\mathbf{x})d\mathbf{%
x}\right) &=&\lim_{m\rightarrow \infty }S^{\prime }\left( \int \theta
(x_{n})dx_{n},S\left( \int \theta (x_{m})dx_{m}\right) \right) \\
&=&\lim_{m\rightarrow \infty }S\left( \int \theta (x_{n})dx_{n}\oplus \int
\theta (x_{m})dx_{m}\right) \\
&=&\lim_{m\rightarrow \infty }S\left( \int \widetilde{\theta }%
(x_{n},x_{m})d(x_{n},x_{m})\right) \\
&=&\int \tilde{\theta}\left( S\left( x_{n}\oplus x_{m}\right) \right)
dS\left( x_{n}\oplus x_{m}\right) \\
&=&\int \widetilde{\theta }(S^{\prime }(x_{n},\mathbf{x}))dS^{\prime }(x_{n},%
\mathbf{x}).
\end{eqnarray*}

In in the last line we have used once again the continuity of the integral
and the assumption that $\left( x_{n}\right) $ is a good $p$-rough path
sequence associated to $\mathbf{x}$.

Therefore, $S^{\prime }\left( \int \theta (x_{n})dx_{n},\int \theta (\mathbf{%
x})d\mathbf{x}\right) $ converges when $n\rightarrow \infty $ to%
\begin{eqnarray*}
\int \widetilde{\theta }(S^{\prime \prime }\left( \mathbf{x})\right)
dS^{\prime \prime }(\mathbf{x}) &=&\lim_{n}\int \widetilde{\theta }%
(S^{\prime \prime }(S\left( x_{n}\right) ))dS^{\prime \prime }(S\left(
x_{n}\right) ) \\
&=&\lim_{n}\int \widetilde{\theta }(S(x_{n}\oplus x_{n}))dS(x_{n}\oplus
x_{n}) \\
&=&\lim_{n}S\left( \int \theta (x_{n})dx_{n}\oplus \int \theta
(x_{n})dx_{n}\right) \\
&=&S^{\prime \prime }\left( \int \theta (\mathbf{x})d\mathbf{x}\right) ,
\end{eqnarray*}%
which proves the first assertion.

For the second one, consider the map%
\begin{equation*}
\left. 
\begin{array}{c}
\widehat{\theta }:\mathbb{R}^{d}\oplus \mathbb{R}^{d}\rightarrow Hom(\mathbb{%
R}^{d}\oplus \mathbb{R}^{d},\mathbb{R}^{n}) \\ 
((x,y),(dX,dY))\rightarrow \theta (y)dX.%
\end{array}%
\right.
\end{equation*}%
\newline
By the continuity of the integral, we obtain that $\int \widehat{\theta }%
(S^{\prime }(x_{n},\mathbf{x}))dS^{\prime }(x_{n},\mathbf{x})$ converges in
the topology induced by $d_{\omega ,p}$ to $\int \widehat{\theta }(S^{\prime
\prime }(\mathbf{x}))dS^{\prime \prime }(\mathbf{x})$. This is our result as 
$\int \widehat{\theta }(S^{\prime }(x_{n},\mathbf{x}))dS^{\prime }(x_{n},%
\mathbf{x})=S\left( \int \theta (\mathbf{x})dx_{n}\right) $.
\end{proof}

We show now that good rough path sequences are stable under the It\^{o} map.

\begin{theorem}
\label{staygoodDE}Let $V_{1},\cdots ,V_{d}$ be $d$ elements of $Lip\left(
p+\varepsilon \right) $-vector fields on $\mathbb{R}^{n},$ and $V=\left(
V_{1},\cdots ,V_{d}\right) $ thought as a linear map from $\mathbb{R}^{d}$
into $Lip\left( p+\varepsilon \right) $ vector fields on $\mathbb{R}^{n}$.
Assume that $\left( x_{n}\right) _{n}$ is a good $p$-rough path sequence
associated to $\mathbf{x}$ for the control $\omega $. Denote by $y_{n}$ the
solution of the ordinary differential equation%
\begin{equation*}
\left\{ 
\begin{array}{l}
dy_{n}(t)=V(y_{n}(t))dx_{n}(t) \\ 
y_{n}(0)=y_{0}.%
\end{array}%
\right.
\end{equation*}%
Then $\left( x_{n}\oplus y_{n}\right) _{n}$ is a good $p$-rough path
sequence associated to $\mathbf{z}=I_{y_{0},V}(\mathbf{x})$ for the control $%
\omega $.
\end{theorem}

\begin{proof}
The proof is similar to the proof of the theorem \ref{staygood}. Denote by $%
\widetilde{V}$ the linear map from $\mathbb{R}^{d}\oplus \mathbb{R}^{d}$
into the $Lip\left( p+\varepsilon \right) $ vector fields on $\mathbb{R}^{n}$
by the formula $\widetilde{V}(y_{1},y_{2})\left( dx_{1},dx_{2}\right)
=(V(y_{1})dx_{1},V(y_{2})dx_{2})$. First notice that 
\begin{equation*}
I_{(y_{0},y_{0}),\widetilde{V}}\left( S^{\prime \prime }(\mathbf{x})\right)
=S^{\prime \prime }(I_{y_{0},V}(\mathbf{x}))
\end{equation*}%
and that 
\begin{equation*}
S^{\prime }(z_{n},\mathbf{z})=I_{(y_{0},y_{0}),\widetilde{V}}(S^{\prime
}(x_{n},\mathbf{x})).
\end{equation*}%
Hence, since $x_{n}$ is a good sequence and by continuity of the It\^{o}
map, $S^{\prime }(z_{n},\mathbf{z})$ converges as $n\rightarrow \infty $ to%
\begin{equation*}
I_{(y_{0},y_{0}),\widetilde{V}}\left( S^{\prime \prime }(\mathbf{x})\right)
=S^{\prime \prime }(I_{y_{0},V}(\mathbf{x}))=S^{\prime \prime }\left( 
\mathbf{z}\right) .
\end{equation*}
\end{proof}

From our two previous theorems, we immediately obtain the following
corollary:

\begin{corollary}
\label{for_Strat}We keep the notation of theorem \ref{staygoodDE}. Then, in
the topology induced by $d_{\omega ,p}$,%
\begin{equation*}
\mathbf{y}=\lim_{n\rightarrow \infty }S\left( y_{0}+\int V\left( \mathbf{y}%
_{u}^{1}\right) dx_{n}(u)\right) .
\end{equation*}%
In particular, looking at the first level of this equation, we obtain that 
\begin{equation*}
y_{0}+\int_{0}^{.}V(y_{u})dx_{n}(u)\rightarrow y.
\end{equation*}
\end{corollary}

\begin{remark}
In the previous theorem and its corollary, with no modification in the
proof, for $2\leq p<3$, one can obtain the same results replacing the map $%
I_{y_{0},V}$ by the map $I_{y_{0},\left( V_{0},V\right) }$, where $V_{0}$ is
a $Lip\left( 1+\varepsilon \right) $-vector fields on $\mathbb{R}^{n}$. In
other words, one can consider differential equations with a time drift and
(almost) minimal smoothness condition on $V_{0}$.
\end{remark}

\subsection{Piecewise-linear Approximation of Brownian motion as a Good
Rough Path Sequence.}

We fix a $p\in (2,3)$ and $\omega (s,t)=t-s$. We recall from the
introduction that $B$ is a $d$-dimensional Brownian motion, and that $%
\mathbf{B}$ is the Stratonovich lift of $B$ to a geometric $p$-rough path.
Let $B^{D}$ be the $D$-linear approximation of $B$ (equation (\ref{Dlinear}%
)). Let $D^{n}$ a sequence of subdivisions of $[0,1]$ which steps tends to $%
0 $ when $n$ tends to infinity. If $D^{n}=\left( \frac{k}{2^{n}},0\leq k\leq
2^{n}\right) ,$ we know from \cite{LQ} that, almost surely, $S\left(
B^{D^{n}}\right) $ converges in $p$-variation to $\mathbf{B}$. If $D^{n}$ is
an increasing sequence of subdivision, i.e. if $D^{n+1}\subset D^{n}$ for
all $n$, a martingale argument proved in \cite{FV1} that almost surely and
in $L^{q}$, $q\geq 1$, $S\left( B^{D^{n}}\right) $ converges in $1/p$-H\"{o}%
lder distance (and even some modulus distances) to $\mathbf{B}$. The
following theorem goes a bit deeper in the convergence of piecewise linear
approximations of the Brownian motion.

\begin{theorem}
\label{Bngood}Let $D^{n}$ be a sequence of subdivision which steps size
tends to $0$. Then $d_{\omega ,p}\left( S^{\prime }\left( B^{D^{n}},\mathbf{B%
}\right) ,S^{\prime \prime }\left( \mathbf{B},\mathbf{B}\right) \right) $
converges when $n$ tends to infinity to $0$ in $L^{q}$, $q\geq 1$ and in
probability.\newline
If $D^{n}=\left( \frac{k}{2^{n}},0\leq k\leq 2^{n}\right) $, the convergence
also holds almost surely, i.e. $B^{n}=B^{D^{n}}$ is almost surely a good $p$%
-rough path sequence associated to $\mathbf{B}$.
\end{theorem}

We decompose the proof in four lemmas.

\begin{lemma}
\label{lemma1}For all $q>p>1,$ the norm $L^{q}$ and $L^{p}$ norm on the $%
k^{th}$ Wiener chaos.
\end{lemma}

\begin{proof}
This is a simple consequence of the hypercontractivity of the
Ornstein-Uhlenbeck semigroup, see \cite[p57]{Nu} for example.
\end{proof}

\begin{lemma}
Let $D$ be a subdivision of $[0,1]$. Then, for all $q\geq 1$ and $p^{\prime
}>1/H,$ there exists $\mu >0$ and $C_{p^{\prime },q,\mu }<\infty $ such that
for all $s<t\in D$%
\begin{equation*}
\left\Vert \int_{s}^{t}B_{s,u}\otimes dB_{u}^{D}-\mathbf{B}%
_{s,t}^{2}\right\Vert _{L^{q}}\leq C_{p^{\prime },q,\mu }\left\vert
D\right\vert ^{\mu }\left\vert t-s\right\vert ^{2/p^{\prime }}
\end{equation*}%
\begin{equation*}
\left\Vert \int_{s}^{t}B_{s,u}^{D}\otimes B_{s,u}^{D}-\mathbf{B}%
_{s,t}^{2}\right\Vert _{L^{q}}\leq C_{p^{\prime },q,\mu }\left\vert
D\right\vert ^{\mu }\left\vert t-s\right\vert ^{2/p^{\prime }}
\end{equation*}
\end{lemma}

\begin{proof}
From the previous lemma, we can take $q=2$. We write $\left( s,t\right)
=\left( t_{m},t_{n}\right) ,$ with $0\leq m<n\leq \left\vert D\right\vert ,$
where $D=\left( t_{i}\right) _{0\leq i\leq D}$. It is easy to see that%
\begin{eqnarray*}
\int_{s}^{t}B_{s,u}\otimes dB_{u}^{D}-\mathbf{B}_{s,t}^{2}
&=&\sum_{k=m}^{n-1}\left( \int_{t_{k}}^{t_{k+1}}B_{t_{k},u}\otimes
dB_{u}^{D}-\mathbf{B}_{t_{k},t_{k+1}}^{2}\right) \\
\int_{s}^{t}B_{s,u}^{D}\otimes B_{s,u}^{D}-\mathbf{B}_{s,t}^{2}
&=&\sum_{k=m}^{n-1}\left( \left( B_{t_{k},t_{k+1}}\right) ^{\otimes 2}-%
\mathbf{B}_{t_{k},t_{k+1}}^{2}\right) .
\end{eqnarray*}%
Therefore, by independence of incremement,%
\begin{eqnarray*}
\left\Vert \int_{s}^{t}B_{s,u}\otimes dB_{u}^{D}-\mathbf{B}%
_{s,t}^{2}\right\Vert _{L^{2}}^{2} &=&\sum_{k=m}^{n-1}E\left( \left(
\int_{t_{k}}^{t_{k+1}}B_{t_{k},u}\otimes dB_{u}^{D}-\mathbf{B}%
_{t_{k},t_{k+1}}^{2}\right) ^{2}\right) \\
&=&C\sum_{k=m}^{n-1}\left( t_{k+1}-t_{k}\right) ^{2} \\
&\leq &C\left\vert D\right\vert ^{2-4/p^{\prime }}\sum_{k=m}^{n-1}\left(
t_{k+1}-t_{k}\right) ^{4/p^{\prime }} \\
&\leq &C\left\vert D\right\vert ^{2-4/p^{\prime }}\left( t-s\right)
^{4/p^{\prime }}.
\end{eqnarray*}%
We also obtain the same estimate for $\left\Vert
\int_{s}^{t}B_{s,u}^{D}\otimes B_{s,u}^{D}-\mathbf{B}_{s,t}^{2}\right\Vert
_{L^{q}},$ which concludes the proof.
\end{proof}

\begin{lemma}
Let $D$ be a subdivision of $[0,1]$. Then, for all $q\geq 1,$ there exists $%
\mu >0$ and $C_{q,\mu }<\infty $ such that%
\begin{equation}
\left\Vert \sup_{0\leq s<t\leq 1}\frac{\left\vert \int_{s}^{t}B_{s,u}\otimes
dB_{u}^{D}-\mathbf{B}_{s,t}^{2}\right\vert }{\left\vert t-s\right\vert ^{2/p}%
}\right\Vert _{L^{q}}\leq C_{q,\mu }\left\vert D\right\vert ^{\mu }
\label{sansD}
\end{equation}%
\begin{equation}
\left\Vert \sup_{0\leq s<t\leq 1}\frac{\left\vert
\int_{s}^{t}B_{s,u}^{D}\otimes B_{s,u}^{D}-\mathbf{B}_{s,t}^{2}\right\vert }{%
\left\vert t-s\right\vert ^{2/p}}\right\Vert _{L^{q}}\leq C_{q,\mu
}\left\vert D\right\vert ^{\mu }  \label{2D}
\end{equation}
\end{lemma}

\begin{proof}
We only prove equation (\ref{sansD}) as the proof for the other estimates is
similar. We define $X_{s,t}^{D}=\int_{s}^{t}B_{s,u}\otimes dB_{u}^{D}-%
\mathbf{B}_{s,t}^{2}$, where $D=\left( 0=t_{0}\leq t_{1}<\cdots
<t_{\left\vert D\right\vert }=1\right) $ is fixed subdivision of $[0,1]$.
First assume that $t_{i}\leq s<t\leq t_{i+1}$. We let $p^{\prime }\in (2,p)$%
; it is easy to check that there exists $C<\infty $ independent of $D$ such
that $\left\Vert B^{D}\right\Vert _{\omega ,p^{\prime }}\leq C\left\Vert
B\right\Vert _{\omega ,p^{\prime }}$. Therefore,%
\begin{eqnarray*}
\frac{\left\vert X_{s,t}^{D}\right\vert }{\left\vert t-s\right\vert ^{2/p}}
&=&\frac{\left\vert \mathbf{B}_{s,t}^{2}\right\vert +\left\vert \left(
\int_{s}^{t}B_{s,u}du\right) \otimes \frac{B_{t_{i},t_{i+1}}}{t_{i+1}-t_{i}}%
\right\vert }{\left\vert t-s\right\vert ^{2/p}} \\
&\leq &\frac{C\left\Vert B\right\Vert _{\omega ,p^{\prime }}^{2}\left(
\left\vert t-s\right\vert ^{2/p^{\prime }}+\left(
\int_{s}^{t}(u-s)^{1/p^{\prime }}du\right) \left( t_{i+1}-t_{i}\right)
^{1/p^{\prime }-1}\right) }{\left\vert t-s\right\vert ^{2/p}} \\
&\leq &C\left\Vert B\right\Vert _{\omega ,p^{\prime }}^{2}\left( \left\vert
t-s\right\vert ^{2/p^{\prime }-2/p}+\left( t-s\right) ^{1+1/p^{\prime
}-2/p}\left( t_{i+1}-t_{i}\right) ^{1/p^{\prime }-1}\right) .
\end{eqnarray*}%
Bounding $t-s$ and $t_{i+1}-t_{i}$ by $\left\vert D\right\vert $, the mesh
size of $D$, we obtain that%
\begin{equation}
\max_{k\in \left\{ 0,\cdots ,\left\vert D\right\vert \right\}
}\sup_{t_{k}\leq s<t\leq t_{k+1}}\frac{\left\vert X_{s,t}^{D}\right\vert }{%
\left\vert t-s\right\vert ^{2/p}}\leq C\left\Vert B\right\Vert _{\omega
,p^{\prime }}^{2}\left\vert D\right\vert ^{2\left( \frac{1}{p^{\prime }}-%
\frac{1}{p}\right) }.  \label{sandtin}
\end{equation}%
Defining $t_{D}$ to be the biggest real in $D$ less than or equal to $t$,
and $s^{D}$ the smallest real in $D$ greater than to $s.$ The above estimate
rewrites%
\begin{equation}
\sup_{\substack{ 0\leq s<t\leq 1  \\ t_{D}<s^{D}}}\frac{\left\vert
X_{s,t}^{D}\right\vert }{\left\vert t-s\right\vert ^{2/p}}\leq C\left\Vert
B\right\Vert _{\omega ,p^{\prime }}^{2}\left\vert D\right\vert ^{2\left( 
\frac{1}{p^{\prime }}-\frac{1}{p}\right) }  \label{sandtin2}
\end{equation}%
and an $L^{q}\,$-estimate is immediate. Hence we are left to prove that%
\begin{equation*}
\left\Vert \sup_{\substack{ 0\leq s<t\leq 1  \\ s^{D}\leq t_{D}}}\frac{%
\left\vert X_{s,t}^{D}\right\vert }{\left\vert t-s\right\vert ^{2/p}}%
\right\Vert _{L^{q}}\leq C_{\mu }\left\vert D\right\vert ^{\mu }.
\end{equation*}%
Now observe that for all $s<t<u$, 
\begin{equation}
X_{s,u}^{D}=X_{s,t}^{D}+X_{t,u}^{D}+B_{s,t}\otimes \left(
B_{t,u}^{D}-B_{t,u}\right) .  \label{badadditivity}
\end{equation}%
Hence, for all $s<t$ such that $s^{D}\leq t_{D}$, so that%
\begin{equation*}
X_{s,t}^{D}=X_{s,s^{D}}^{D}+X_{s^{D},t_{D}}^{D}+X_{t_{D},t}^{D}+B_{s,t_{D}}%
\otimes \left( B_{t_{D},t}^{D}-B_{t_{D},t}\right) .
\end{equation*}

From (\ref{sandtin2})) 
\begin{equation*}
\left\Vert \sup_{\substack{ 0\leq s<t\leq 1  \\ s^{D}\leq t_{D}}}\frac{%
\left\vert X_{s,s^{D}}^{D}\right\vert }{\left\vert s^{D}-s\right\vert ^{2/p}}%
+\frac{\left\vert X_{t_{D},t}^{D}\right\vert }{\left\vert t-t_{D}\right\vert
^{2/p}}\right\Vert _{L^{q}}\leq C_{\mu }\left\vert D\right\vert ^{\mu }.
\end{equation*}%
Compatibility of tensor norms shows - similar as above but easier -%
\begin{equation*}
\frac{\left\vert B_{s,t_{D}}\otimes \left(
B_{t_{D},t}^{D}-B_{t_{D},t}\right) \right\vert }{\left\vert t-s\right\vert
^{2/p}}\leq C\left\Vert B\right\Vert _{\omega ,p^{\prime }}^{2}\left\vert
D\right\vert ^{1/p^{\prime }-1/p};
\end{equation*}%
therefore, we just need to check that for some $\mu >0$,%
\begin{equation}
\left\Vert \sup_{0\leq s<t\leq 1}\frac{\left\vert
X_{s^{D},t_{D}}^{D}\right\vert }{\left\vert t-s\right\vert ^{2/p}}%
\right\Vert _{L^{q}}\overset{\text{trivial}}{\leq }\left\Vert \max_{s<t\in D}%
\frac{\left\vert X_{s,t}^{D}\right\vert }{\left\vert t-s\right\vert ^{2/p}}%
\right\Vert _{L^{q}}\leq C_{\mu }\left\vert D\right\vert ^{\mu }.
\label{whattoprove}
\end{equation}%
Now consider $\widetilde{X}^{D}:\left[ 0,1\right] \rightarrow \mathbb{R}%
^{d}\otimes \mathbb{R}^{d}$ the linear path on the intervals $\left[
t_{i},t_{i+1}\right] $, such that for all $i$, $\widetilde{X}%
_{t_{i},t_{i+1}}^{D}=X_{t_{i},t_{i+1}}^{D}$, i.e. for all $t,$ 
\begin{equation*}
\widetilde{X}_{t_{D},t}^{D}=X_{t_{i},t}^{D}+\frac{t-t_{D}}{t^{D}-t_{D}}%
X_{t_{D,}t^{D}}^{D}.
\end{equation*}%
The previous lemma showed that for all $s,t\in D,$ $\left\Vert \widetilde{X}%
_{s,t}^{D}\right\Vert _{L^{q}}\leq C_{p^{\prime },q,\mu }\left\vert
D\right\vert ^{\mu }\left\vert t-s\right\vert ^{2/p^{\prime }}.$ From this,
it is easy to check that (changing the constant), this equality remains true
for all $s,t\in \left[ 0,1\right] $. Define $B_{\alpha
}=\int_{0}^{1}\int_{0}^{1}\left\vert \frac{\widetilde{X}_{s,t}^{D}}{\left(
t-s\right) ^{\alpha }}\right\vert ^{2m}dsdt$. If $\alpha <\frac{2}{p^{\prime
}}+\frac{1}{2m}$, then $\mathbb{E}(B_{\alpha })\leq C_{m}\left\vert
D\right\vert ^{\mu m}$. By Garsia, Rodemich and Rumsey's theorem \cite{St},
we obtain that%
\begin{equation*}
\sup_{0\leq s<t\leq 1}\frac{\left\vert \widetilde{X}_{s,t}^{D}\right\vert }{%
\left( t-s\right) ^{\alpha -1/m}}\leq CB_{\alpha }^{1/2m}.
\end{equation*}%
Therefore, there exists two constants $\vartheta >0$ and $C_{\vartheta
}<\infty $ (independent of $D$) such that%
\begin{equation*}
\left\Vert \max_{s<t\in D}\frac{\left\vert X_{s,t}^{D}\right\vert }{%
\left\vert t-s\right\vert ^{2/p}}\right\Vert _{L^{q}}\leq C_{\vartheta
}\left\vert D\right\vert ^{\vartheta }.
\end{equation*}
\end{proof}

\begin{lemma}
Let $D$ be a subdivision of $[0,1]$. Then, for all $q\geq 1,$ there exists $%
\nu >0$ and $C_{q,\nu }<\infty $ such that%
\begin{equation*}
\left\Vert d_{\omega ,p}\left( S^{\prime }\left( B^{D},\mathbf{B}\right)
,S^{\prime \prime }\left( \mathbf{B},\mathbf{B}\right) \right) \right\Vert
_{L^{q}}\leq C_{q,\nu }\left\vert D\right\vert ^{\nu }.
\end{equation*}
\end{lemma}

\begin{proof}
It is easy to check that $\sup_{0\leq s<t\leq 1}\frac{\left\vert \mathbf{B}%
_{s,t}^{1}-B_{s,t}^{D}\right\vert }{\left\vert t-s\right\vert ^{1/p}}\leq
C\left\Vert \mathbf{B}\right\Vert _{\omega ,q}\left\vert D\right\vert
^{1/q-1/p},$where $2<q<p$. The previous lemmae together with proposition \ref%
{good2} give the result.
\end{proof}

We can now turn to the proof of theorem \ref{Bngood}.

\begin{proof}
The first part is obvious from the previous lemma. From the second part,%
\begin{eqnarray*}
\mathbb{P}\left( d_{\omega ,p}\left( S^{\prime }\left( B^{n},\mathbf{B}%
\right) ,S^{\prime \prime }\left( \mathbf{B},\mathbf{B}\right) \right) \geq 
\frac{1}{n}\right) &\leq &n^{2}\left\Vert d_{\omega ,p}\left( S^{\prime
}\left( B^{n},\mathbf{B}\right) ,S^{\prime \prime }\left( \mathbf{B},\mathbf{%
B}\right) \right) \right\Vert _{L^{2}}^{2} \\
&\leq &C_{q,\nu }n^{2}2^{-n\nu }.
\end{eqnarray*}%
Hence, by Borel-Cantelli's lemma, we obtain that, almost surely, $S\prime
\left( B^{n},\mathbf{B}\right) $ converges to $S^{\prime \prime }\left( 
\mathbf{B},\mathbf{B}\right) $ in the topology induced by $d_{\omega ,p}$.
\end{proof}

\subsection{\protect\bigskip Piecewise-linear Approximation of Fractional
Brownian motion as a Good Rough Path Sequence.}

Let $W_{H}$ be fractional Browian motion of Hurst paramter $H\in (0,1),~$%
defined as Gaussian process on $[0,1]$ with zero mean and covariance%
\begin{equation*}
\mathbb{E}\left( W_{H}\left( t\right) W_{H}\left( s\right) \right) =\frac{1}{%
2}\left[ t^{2H}+s^{2H}-\left\vert t-s\right\vert ^{2H}\right] .
\end{equation*}%
From Kolmogorov's criterion, $W_{H}$ has "$H-\epsilon $"-H\"{o}lder
regularity for any $\epsilon \in (0,H).$ For $H>1/2$ this is enough
regularity to determine all iterated integrals as Young-integrals. Then a
pathwise ODE-theory is possible, see \cite{LQ}. When $H=1/2$ we are dealing
with usual Brownian motion. For $H\in \left( \frac{1}{4},\frac{1}{2}\right) $
the canonical lift to a geometric rough path $\mathbf{W}_{H}$ was
constructed in \cite{CQ}. According to rough path theory $H\in (\frac{1}{3},%
\frac{1}{2})$ requires $\mathbf{W}_{H}=\left( 1,\mathbf{W}_{H}^{1},\mathbf{W}%
_{H}^{2}\right) $ while for $H\in (\frac{1}{4},\frac{1}{3}]$ we need
additionally the third iterated integrals,%
\begin{equation*}
\mathbf{W}_{H}=\left( 1,\mathbf{W}_{H}^{1},\mathbf{W}_{H}^{2},\mathbf{W}%
_{H}^{3}\right) .
\end{equation*}%
As in the Brownian case, we work with the H\"{o}lder control $\omega \left(
s,t\right) =t-s$ and take%
\begin{equation*}
1/H<p<\lfloor 1/H\rfloor +1.
\end{equation*}%
The proofs of the following two theorems are found in the appendix.

\begin{theorem}
\label{CQ copy(1)}In $d_{\omega ,p}$-topology%
\begin{equation*}
S\left( W^{D^{n}}\right) \rightarrow \mathbf{W}_{H}
\end{equation*}
in $L^{q}$ for any $q\in \lbrack 1,\infty )$ and in probability.
\end{theorem}

\begin{theorem}
\label{Bngoodmbf copy(1)}In $d_{\omega ,p}$-topology%
\begin{equation*}
S^{\prime }\left( W_{H}^{D^{n}},\mathbf{W}\right) \rightarrow S^{\prime
\prime }\left( \mathbf{W}_{H},\mathbf{W}_{H}\right)
\end{equation*}%
in $L^{q}$ for any $q\in \lbrack 1,\infty )$ and in probability. If $%
D^{n}=\left( \frac{k}{2^{n}},0\leq k\leq 2^{n}\right) $, the convergence
also holds almost surely. In other words, $\left( W_{H}^{D^{n}}\right)
_{n\in \mathbb{N}}$ is a.s. a good $p$-rough path sequence associated to $%
\mathbf{W}_{H}$.\bigskip\ 
\end{theorem}

\section{Anticipative and Fractional Stochastic Analysis}

Combing the results of the preceding sections yields a powerful theory which
unifies anticipation with non-semimartingale driving noise such as
fractional Brownian motion. For the reader's convenience, we always discuss
the case of driving Brownian motion separately (although the arguments we
developed allow for rather similar proofs). As concrete applications, we
discuss Wong-Zakai results, support theorems and large deviations in the
simultaneous context of anticipation and driving fractional Brownian motion.

\subsection{Rough Paths Solution Equals Stratonovich Solution}

\subsubsection{The case of driving Brownian motion}

\bigskip Fix $p\in (2,3)$ and $\omega (s,t)=t-s$. Consider random
vector-fields $V_{0}$ and $V=\left( V_{1},\cdots ,V_{d}\right) ,$ where $%
V_{0}$ (resp. $V_{1},\cdots ,V_{d}$) is almost surely a $Lip\left(
1+\varepsilon \right) $ (resp. $Lip\left( 2+\varepsilon \right) $) vector
field on $\mathbb{R}^{n}$, and an a.s. finite random variable $y_{0}\in 
\mathbb{R}^{n}$. As earlier, $\mathbf{B}$ denotes the canonical lift of
Brownian motion to a geometric $p$-rough path. Then, there exists a unique
rough path solution of the (anticipative) stochastic differential equation%
\begin{equation*}
d\mathbf{y}_{t}\mathbf{=}V_{0}\left( \mathbf{y}_{t}\right) dt+V\left( 
\mathbf{y}_{t}\right) d\mathbf{B}_{t}
\end{equation*}%
with (random) initial condition $y_{0}$.

\begin{theorem}
Let $\mathbf{y}^{1}$ denote the projection to path-level of the rough-path $%
\mathbf{y.}$ Then $\mathbf{y}^{1}$ solves the anticipative Stratonovich
stochastic differential equation,%
\begin{equation}
\mathbf{y}_{t}^{1}=y_{0}+\int_{0}^{t}V_{0}(\mathbf{y}_{u}^{1})du+%
\int_{0}^{t}V(\mathbf{y}_{u}^{1})\circ dB_{u}.
\end{equation}%
Moreover, almost surely, $\mathbf{y}^{1}$ is the limit in $1/p$-H\"{o}lder
topology of 
\begin{equation*}
t\rightarrow y_{0}+\int_{0}^{t}V_{0}(\mathbf{y}_{u}^{1})du+\int_{0}^{t}V(%
\mathbf{y}_{u}^{1})dB_{u}^{n},
\end{equation*}%
where $B^{n}$ is the dyadic linear approximation of $B$ of level $n$.
\end{theorem}

\begin{proof}
For any sequence of subdivisions $\left( D^{n}\right) _{n}$ which mesh size
tends to $0$ when $n\rightarrow \infty $, $d_{\omega ,p}\left( S^{\prime
}\left( B^{D^{n}},\mathbf{B}\right) ,S^{\prime \prime }\left( \mathbf{B},%
\mathbf{B}\right) \right) $ converges in probability to $0$ (theorem \ref%
{Bngood}) hence by corollary \ref{for_Strat},%
\begin{equation*}
d_{\omega ,p}\left( S\left( y_{0}+\int_{0}^{t}V_{0}(\mathbf{y}%
_{u}^{1})du+\int_{0}^{t}V(\mathbf{y}_{u}^{1})dB_{u}^{D^{n}}\right) ,\mathbf{y%
}\right)
\end{equation*}%
converges in probability to $0$. The first level of this equation says
precisely that $\mathbf{y}^{1}$ is solution of the anticipative Stratonovich
stochastic differential equation (\ref{ASDE}).

The same argument and the fact that $d_{\omega ,p}\left( S^{\prime }\left(
B^{n},\mathbf{B}\right) ,S^{\prime \prime }\left( \mathbf{B},\mathbf{B}%
\right) \right) $ converges almost surely to $0$ gives the second part of
the theorem.
\end{proof}

\subsubsection{The case of driving fBM with $H\in (1/4,1/2]$}

\bigskip Of course, this section covers the one above with $H=1/2.$We fix $%
H\in (1/4,1/2]$ and $p\in (1/H,\lfloor 1/H\rfloor +1)$. As earlier, $\mathbf{%
W}_{H}$ denotes the $\lfloor p\rfloor $-level rough path associated to
fractional Brownian Motion and we work with H\"{o}lder-control $\omega
(s,t)=t-s.$

Consider random vectorfields $V_{0}$ and $V=\left( V_{1},\cdots
,V_{d}\right) ,$ where $V_{0}$ (resp. $V_{1},\cdots ,V_{d}$) is almost
surely a $Lip\left( 1+\varepsilon \right) $ (resp. $Lip\left(
1/H+\varepsilon \right) $) vector field on $\mathbb{R}^{n}$, and an a.s.
finite random variable $y_{0}\in \mathbb{R}^{n}$. Then, there exists a
unique rough path solution of the (anticipative) stochastic differential
equation%
\begin{equation*}
d\mathbf{y}_{t}\mathbf{=}V_{0}\left( \mathbf{y}_{t}\right) dt+V\left( 
\mathbf{y}_{t}\right) d\mathbf{W}_{H}\left( t\right)
\end{equation*}%
with (random) initial condition $y_{0}$.

\begin{theorem}
$\mathbf{y}^{1}$ solves the anticipative Stratonovich stochastic
differential equation%
\begin{equation}
\mathbf{y}_{t}^{1}=y_{0}+\int_{0}^{t}V_{0}(\mathbf{y}_{u}^{1})du+%
\int_{0}^{t}V(\mathbf{y}_{u}^{1})\circ dW_{H}\left( u\right) .  \label{ASDE}
\end{equation}%
Moreover, almost surely, $\mathbf{y}^{1}$ is the limit in $1/p$-H\"{o}lder
topology of 
\begin{equation*}
t\rightarrow y_{0}+\int_{0}^{t}V_{0}(\mathbf{y}_{u}^{1})du+\int_{0}^{t}V(%
\mathbf{y}_{u}^{1})dW_{H}^{n}\left( u\right) ,
\end{equation*}%
where $W_{H}^{n}$ is the $n^{th}$ dyadic linear approximation of $W_{H}$.
\end{theorem}

\begin{proof}
For any sequence of subdivisions $\left( D^{n}\right) _{n}$ which mesh size
tends to $0$ when $n\rightarrow \infty $, $d_{\omega ,p}\left( S^{\prime
}\left( W_{H}^{D^{n}},\mathbf{W}\right) ,S^{\prime \prime }\left( \mathbf{W}%
_{H},\mathbf{W}_{H}\right) \right) $ converges in probability to $0$ by
theorem \ref{Bngoodmbf copy(1)} hence by corollary \ref{for_Strat},%
\begin{equation*}
d_{\omega ,p}\left( S\left( y_{0}+\int_{0}^{t}V_{0}(\mathbf{y}%
_{u}^{1})du+\int_{0}^{t}V(\mathbf{y}_{u}^{1})dB_{u}^{D^{n}}\right) ,\mathbf{y%
}\right)
\end{equation*}%
converges in probability to $0$. The first level of this equation says
precisely that $\mathbf{y}^{1}$ is solution of the anticipative Stratonovich
stochastic differential equation (\ref{ASDE}).

The same argument and the fact that $d_{\omega ,p}\left( S^{\prime }\left(
W_{H}^{^{n}},\mathbf{W}\right) ,S^{\prime \prime }\left( \mathbf{W}_{H},%
\mathbf{W}_{H}\right) \right) $ converges almost surely to $0$ gives the
second part of the theorem.
\end{proof}

\subsection{A Wong-Zakai Theorem}

The Universal Limit Theorem gives us for free a\ Wong-Zakai theorem for our
solution of the Stratonovich differential equation. The situation is so
simple that we feel no need to discuss $H=1/2$ seperately.

\begin{theorem}
Under the same assumptions and notation than the previous theorem, let $%
y^{n} $ be the solution of the differential equation%
\begin{equation*}
\left\{ 
\begin{array}{l}
dy_{t}^{n}=V_{0}(y_{t}^{n})dt+V\left( y_{t}^{n}\right) dW_{H}^{n}\left(
t\right) , \\ 
y_{0}^{n}=y_{0}.%
\end{array}%
\right.
\end{equation*}%
Then almost surely, $y^{n}$ converges to $y$ in the topology induced by the $%
1/p$-H\"{o}lder distance.
\end{theorem}

\begin{proof}
Let $\mathbf{z}^{n}=I_{V,y_{0}}(S(W_{H}^{n}\mathbf{))})$. Then by the
continuity of the It\^{o} map, $d_{\omega ,p}(\mathbf{z}^{n},\mathbf{z}%
)\rightarrow 0$. As $\mathbf{z}^{n}$ projects down onto $y^{n}$, this is a
stronger result that the stated theorem.
\end{proof}

\subsection{A\ Large Deviation Principle}

\subsubsection{The case of driving Brownian motion}

Due to the universal limit theorem, Freidlin-Wentzell's type theorems have
easy proofs via rough paths \cite{LQZ,FV}. We give here an extension of the
Freidlin-Wentzell's theorem, and of the main theorem in \cite{MNS}.

\begin{theorem}
Let $\left( y_{0}^{\alpha }\right) _{\alpha \geq 0}$ be a family of random
elements of $\mathbb{R}^{d}$, $\left( V_{0}^{\alpha }\right) _{\alpha \geq
0} $ be a family of random $Lip\left( 1+\varepsilon \right) $-vector fields, 
$\left( V_{1}^{\alpha },\cdots ,V_{d}^{\alpha }\right) _{\alpha \geq 0}$, $%
i=1,\cdots ,d$ d families of random $Lip\left( 2+\varepsilon \right) $%
-vector fields, such that for all $\beta >0$ 
\begin{eqnarray*}
\lim_{\alpha \rightarrow 0}\alpha \log \mathbb{P}\left( \left\Vert
V_{0}^{\alpha }-V_{0}^{0}\right\Vert _{Lip\left( 1+\varepsilon \right)
}>\beta \right) &=&-\infty , \\
\lim_{\alpha \rightarrow 0}\alpha \log \mathbb{P}\left( \max_{1\leq i\leq
d}\left\Vert V_{i}^{\alpha }-V_{i}^{0}\right\Vert _{Lip\left( 2+\varepsilon
\right) }>\beta \right) &=&-\infty , \\
\lim_{\alpha \rightarrow 0}\alpha \log \mathbb{P}\left( \left\vert
y_{0}^{\alpha }-y_{0}^{0}\right\vert >\beta \right) &=&-\infty .
\end{eqnarray*}%
Define $\mathbf{y}^{\alpha }$ to be the rough path solution of the
differential equation%
\begin{equation*}
\left\{ 
\begin{array}{l}
d\mathbf{y}_{t}^{\alpha }=V_{0}^{\alpha }\left( \mathbf{y}_{t}^{\alpha
}\right) dt+\sqrt{\alpha }V^{\alpha }\left( \mathbf{y}_{t}^{\alpha }\right) d%
\mathbf{B}_{t} \\ 
\mathbf{y}_{0}^{\alpha }=y_{0}^{\alpha }.%
\end{array}%
\right.
\end{equation*}%
Then $\left( \mathbf{y}^{\alpha }\right) _{\alpha >0}$ satisfies a large
deviation principle in the topology induced by $d_{\omega ,p}$ with good
rate function 
\begin{equation*}
J(\mathbf{x)}=\inf_{I_{y_{0}^{0},V^{0}}\mathbf{(y)=x}}I(\mathbf{y}),
\end{equation*}%
where%
\begin{equation*}
I(\mathbf{x)=}\left\{ 
\begin{array}{l}
\frac{1}{2}\int_{0}^{1}\left\vert x_{u}^{\prime }\right\vert ^{2}du\text{,
if }S(x)=\mathbf{x}\text{ for some }x\in W^{1,2} \\ 
+\infty \text{ otherwise.}%
\end{array}%
\right.
\end{equation*}
\end{theorem}

\begin{proof}
In \cite{FV1}, it was proved that $\left( \delta _{\sqrt{\alpha }}\mathbf{B}%
\right) _{\alpha >0}$ satisfies a large deviation principle in the topology
induced by $d_{\omega ,p}$ with good rate function $I$. The assumptions on
the vector fields and the initial conditions give that $\left( \left(
V_{i}^{\alpha }\right) _{0\leq i\leq d},y_{0}^{\alpha },\delta _{\sqrt{%
\alpha }}\mathbf{B}\right) $ satisfies a large deviation (in the topology
product $Lip\left( 1+\varepsilon \right) ,Lip\left( 2+\varepsilon \right)
^{d},\left\vert .\right\vert ,$ and $d_{\omega ,p}$) with good rate function 
\begin{equation*}
\left\{ 
\begin{array}{l}
\frac{1}{2}\int_{0}^{1}\left\vert x_{u}^{\prime }\right\vert ^{2}du\text{,
if }S(x)=\mathbf{x}\text{ for some }x\in W^{1,2},\text{ }V_{i}=V_{i}^{0}%
\text{, }1\leq i\leq d\text{ and }y_{0}=y_{0}^{0} \\ 
+\infty \text{ otherwise.}%
\end{array}%
\right.
\end{equation*}%
By the continuity of the It\^{o} map (theorem \ref{ULT copy(1)},\ref{ITOIC
copy(1)} and \ref{ITOVF copy(1)}), we obtain our large deviation principle.
\end{proof}

\subsubsection{\protect\bigskip The case of driving fBM with $H\in (1/4,1/2]$%
}

Recent work \cite{MS} establishes a large deviation principle for fractional
Brownian motion lifted to a rough path w.r.t. $p$-variation topology with
good rate function%
\begin{equation*}
I_{H}\left( \mathbf{x}\right) =\left\{ 
\begin{array}{c}
\frac{1}{2}\left\Vert i^{-1}\left( x\right) \right\Vert _{\mathcal{H}%
_{H}}^{2}\text{ \ if }S\left( x\right) =\mathbf{x}\text{ for some }x\in
i\left( \mathcal{H}_{H}\right) \\ 
+\infty \text{ otherwise \ \ \ \ \ \ \ \ \ \ \ \ \ \ \ \ \ \ \ \ \ \ \ \ \ \
\ \ \ \ \ \ \ \ \ \ \ \ \ \ \ \ \ \ \ \ \ \ \ \ \ \ }%
\end{array}%
\right.
\end{equation*}%
This shows that the above LD principle is true in $p$-variation topology
with rate function $I_{H}$ when (enhanced)\ Brownian motion $\mathbf{B}$ is
replaced by (enhanced) fractional Brownian motion $\mathbf{W}_{H}.$

\begin{remark}
We are certain the LD\ principle for SDEs driven by (enhanced)\ fractional
Brownian motion holds in H\"{o}lder topology. We shall discuss this in
forthcoming work.
\end{remark}

\subsection{Support theorem}

It turns out that the notion of good rough path sequence leads to neat
support description of fractional Brownian motion lifted to a rough path. In
fact, given that piecewise linear approximations form a good rough path
sequence to the fractional Brownian rough path, our proof is no more
complicated than the well-known proof on the characterization of the support
of the law of the Brownian motion as seen in Revuz and Yor \cite{RY}. The
method here works for any $H\in (1/4,1)$ but one should keep in mind that
for $H>1/2$ no lift is necessary and the support description is trivial. The
case of Brownian motion, $H=1/2,$ has been dealt with several times,
including by the two last named authors of this article (\cite{F}, \cite{FV}%
). Hence, the real interest of the result below is for $H<1/2$. It was also
discussed, with different methods, in recent work \cite{FV1} and \cite{FP}.

\begin{lemma}
The map 
\begin{eqnarray*}
Minus &:&\mathbb{R}^{d}\oplus \mathbb{R}^{d}\rightarrow \mathbb{R}^{d} \\
x,y &\rightarrow &y-x
\end{eqnarray*}%
induces a homeomorphism $Minus$ from $G^{[p]}\left( \mathbb{R}^{d}\oplus 
\mathbb{R}^{d}\right) $ onto $G^{[p]}\left( \mathbb{R}^{d}\right) $.
\end{lemma}

\begin{example}
Let $\left[ p\right] =2.$ Given $Z\in G^{2}\left( \mathbb{R}^{d}\oplus 
\mathbb{R}^{d}\right) $%
\begin{equation*}
Z=\left( 1,\,\,Z^{1;1}\oplus Z^{1;2},\,\,Z^{2;1,1}\oplus Z^{2;1,2}\oplus
Z^{2;2,1}\oplus Z^{2;2,2}\right) .
\end{equation*}%
with $Z^{1;i}\in \mathbb{R}^{d}$, $Z^{2;i,j}\in \mathbb{R}^{d}\otimes 
\mathbb{R}^{d}$. Then 
\begin{equation*}
Minus\left( Z\right) =\left(
1,\,\,Z^{1;2}-Z^{1;1},\,\,Z^{2;2,2}-Z^{2;1,2}-Z^{2;2,1}+Z^{2;1,1}\right) \in
G^{[p]}\left( \mathbb{R}^{d}\right) .
\end{equation*}
\end{example}

The map $Minus$ induces a continuous surjection, still denoted $Minus$, from 
$C^{0,p-var}\left( G^{[p]}(\mathbb{R}^{d}\oplus \mathbb{R}^{d})\right) $
onto $C^{0,p-var}\left( G^{[p]}(\mathbb{R}^{d})\right) $. We define for $%
\left( h,\mathbf{x}\right) $ in the space $C^{0,q-var}\left( \mathbb{R}%
^{d}\right) \times C^{0,p-var}\left( G^{[p]}(\mathbb{R}^{d})\right) $, with $%
q^{-1}+p^{-1}>1$ the translation of $\mathbf{x}$ by $-h$%
\begin{equation*}
T_{-h}\left( \mathbf{x}\right) =Minus\left( S^{^{\prime }}\left( h,\mathbf{x}%
\right) \right) .
\end{equation*}
The following lemma was obtained in \cite{FV1} for $\mathbf{W}_{1/2}=$ $%
\mathbf{B.}$ The generalization is easy.

\begin{lemma}
\label{CMfbm}Let $h$ be a bounded variation path in the Cameron-Martin space
associated to fBM with $H>1/4.$ Then the law of $T_{h}\left( \mathbf{W}%
_{H}\right) $ is equivalent to the law of $\mathbf{W}_{H}$.\newline
\end{lemma}

\begin{theorem}
\label{supB}Let $H\in (1/4,1)$ and $p\in (1/H,\lfloor 1/H\rfloor +1)$. The
support of the law of $\mathbf{W}_{H}$ in the $d_{\omega ,p}$-topology is
the set of path starting at $0$ in $C^{0,\omega ,p}\left( G^{[p]}\left( 
\mathbb{R}^{d}\right) \right) $ where $\omega (s,t)=t-s.$
\end{theorem}

\begin{proof}
Almost surely, $W_{H}\in C^{0,\omega ,p}\left( G^{[p]}\left( \mathbb{R}%
^{d}\right) \right) $ and $W_{H}\left( 0\right) =0$ therefore the support of
its law is included in the set of path starting at $0$ in $C^{0,\omega
,p}\left( G^{[p]}\left( \mathbb{R}^{d}\right) \right) $. Reciprocally, the
support of the law of $\mathbf{W}_{H}$ contains at least one point $\mathbf{x%
}$ such that its sequence of dyadic linear approximation of level $n$ $%
\left( x^{n}\right) _{n}$ is a good $p$ rough path sequence associated to $%
\mathbf{x}$ (due to theorem \ref{Bngoodmbf copy(1)}). By lemma \ref{CMfbm},
the support of the law of $\mathbf{W}_{H}$ contains the $d_{\omega ,p}$%
-closure of $\left\{ T_{h}\left( \mathbf{x}\right) \text{, }h\text{ in the
Cameron-Martin space}\right\} $. As $\left( x^{n}\right) $ is a good $p$%
-rough path sequence associated to $\mathbf{x}$, $T_{-x^{n}}\left( \mathbf{x}%
\right) =Minus\left( S^{\prime }\left( x^{n},\mathbf{x}\right) \right) $
converges to $0$ (in the $d_{\omega ,p}$-topology). From section \ref{plaacm}%
, piecewise linear approximation are Cameron Martin paths and since the
support of a measure is always closed (by definition) it follows that $0$
belongs to the support of the law of $\mathbf{W}_{H}$. Clearly, by lemma \ref%
{CMfbm}, the support contains the closure of the translation of all smooth
paths. Therefore, the support of the law of $\mathbf{W}_{H}$ contains the
closure in the $d_{\omega ,p}$-topology of the set of smooth paths starting
at $0$. This concludes the proof with the results in \cite{FV1}.
\end{proof}

Denote by $I$ the map which maps $\left( y_{0},\left( V_{0},V\right) ,%
\mathbf{x}\right) $ to $I_{y_{0},\left( V_{0},V\right) }\left( \mathbf{x}%
\right) $. The following proposition is an obvious corollary of the the
continuity of the map $I$ (theorems \ref{ULT copy(1)},\ref{ITOIC copy(1)},%
\ref{ITOVF copy(1)}) and of theorem \ref{supB}.

\begin{proposition}
Let $\mathbf{y}$ be the solution of the rough differential equation 
\begin{equation*}
\left\{ 
\begin{array}{l}
d\mathbf{y}_{t}=V_{0}\left( \mathbf{y}_{t}\right) dt+V\left( \mathbf{y}%
_{t}\right) d\mathbf{W}_{H}\left( t\right) \\ 
\mathbf{y}_{0}^{1}=y_{0},%
\end{array}%
\right.
\end{equation*}%
where $V_{0}$ is almost surely a $Lip\left( 1+\varepsilon \right) $ vector
field and $V_{i},$ $i\in \left\{ 1,\cdots ,d\right\} $ are almost surely $%
Lip\left( 1/H+\varepsilon \right) $ vector fields$,$ $y_{0}\in \mathbb{R}%
^{d} $ is almost surely finite. The support of the law of $\mathbf{y}$ in
the $d_{\omega ,p}$-topology is equal to\ the image by the map $I$ of the
support of the law of $\left( y_{0},V,\mathbf{W}_{H}\right) $, in the
product of the Euclidean, Lipschitz and $d_{\omega ,p}$ topology. \newline
\end{proposition}

In particular, if $y_{0}$ and the vector fields $V_{i}$ are deterministic,
in the $d_{\omega ,p}$-topology, the support of the law of $\mathbf{y}$ is
equal to the set $I_{y_{0},\left( V_{0},V\right) }\left( C^{0,\omega
,p}\left( G^{[p]}\left( \mathbb{R}^{d}\right) \right) \right) $, which, at
the first level and specialized to $H=1/2$ is the classical support theorem
of Stroock-Varadhan \cite{SV}.

If $y_{0}$ and the $V_{i}$s are the image by a continuous function of $%
\mathbf{B}$, then the support is still trivially characterized. One could
then ask for more specific conditions on $y_{0}$ and the $V_{i}$s, in the
spirit of \cite{MN}, and obtain a detailed support theorem. Thanks to the
Universal Limit, one would obtain stronger results than in \cite{MN}
(stronger topology and without the assumption of deterministic vector
fields) but we shall not pursue this here.

\bigskip

\section{\protect\bigskip Appendix 1}

\subsection{Linear Approximation of the Fractional Brownian motion as Good
Rough Path Sequences}

We fix $\omega (s,t)=t-s.$

\subsubsection{Fractional Brownian motion framework}

We use the framework of \cite{CCM} or \cite{CQ}. The starting point of the
approach develloped in \cite{CCM} is the following representation of
fractional Brownian motion given by Decreusefond-\"{U}st\"{u}nel, \cite{Due}%
: ${\mathbf{P}}$ almost every where 
\begin{equation}
W_{H}(t)=\int_{0}^{t}K_{H}(t,s)dB_{s},~~~~\forall t\in \lbrack 0,1]
\label{WHviaKernel}
\end{equation}%
where $K_{H}(t,s)=\frac{(t-s)^{H-\frac{1}{2}}}{\Gamma (H+\frac{1}{2})}F(H-%
\frac{1}{2},\frac{1}{2}-H,H+\frac{1}{2},1-\frac{t}{s}),~~~s<t.$\newline
$F$ denotes Gauss hypergeometric function, \cite{Leb} and $\left(
B_{t}\right) _{t\in \lbrack 0,1]}$ is a Brownian motion. According to Lemma
2.7 of \cite{CCM}, the function $t\mapsto K_{H}(t,s)$ is differentiable on $%
]s,+\infty \lbrack ,$ with derivative 
\begin{equation*}
\partial _{t}K_{H}(t,s)=\frac{(s/t)^{H-\frac{1}{2}}}{\Gamma (\frac{1}{2}-H)}%
(s-t)^{H-\frac{1}{2}},~~~0<s<t.
\end{equation*}%
The parameter $H$ will be fixed in $\left( \frac{1}{4},\frac{1}{2}\right) $
(if $H>1/2,$ things are trivial, $H=\frac{1}{2}$ is the Brownian motion
case, while our techniques do not allow us to deal with the case $H\leq
1/4). $

\subsubsection{"Canonical" lift of fractional Brownian motion to a geometric 
$p$-rough path}

Let be $d\in {\mathbf{N}}^{\ast }$ and $B=(B^{1},\cdots ,B^{d})$ be a $d$
dimensional Brownian motion, $H\in ]\frac{1}{4},1[,$ and $%
W_{H}=(W_{H}^{1},\cdots ,W_{H}^{d})$ given by 
\begin{equation*}
W_{H}^{i}(t)=\int_{0}^{t}K_{H}(t,s)dB_{s}^{i},~~~~\forall t\in \lbrack 0,1].
\end{equation*}%
Indeed, $W_{H}$ is a $d$ dimensional fractional Brownian motion with Hurst
parameter $H.$ We fix $p\in (\frac{1}{H},\left[ \frac{1}{H}\right] +1).$

\begin{definition}
For $H>\frac{1}{4},$ we define, according to \cite{CQ}, the lift of $W_{H}$
to a $p$-geometric rough path ${\mathbf{W}}_{H}$ is formally given by%
\begin{align*}
\mathbf{W}_{H}^{1}(0,t)& =W_{H}(t), \\
\mathbf{W}_{H}^{2}(0,t)& =\int_{0}^{t}W_{H}(s)\otimes \circ dW_{H}(s), \\
\mathbf{W}_{H}^{3}(0,t)& =\int_{0}^{t}\mathbf{W}_{H}^{2}(0,s)\otimes \circ
dW_{H}(s).
\end{align*}%
The 2nd level is rigourosly given as%
\begin{eqnarray*}
\mathbf{W}_{H}^{2}(0,t)^{i,j}
&=&\int_{0}^{1}I_{0,t}^{K_{H}}(W_{H}^{i})(u)dB^{j}(u)\text{ for }i\neq j \\
\mathbf{W}_{H}^{2}(0,t)^{i,i} &=&\frac{1}{2}W_{H}^{i}(t)^{2}
\end{eqnarray*}%
\newline
with%
\begin{equation*}
I_{t,-}^{K_{H}}(f)(s)=K_{H}(t,s)f(s)+\int_{s}^{t}(f(u)-f(s))\partial
_{u}K_{H}(u,s)du.
\end{equation*}%
From \cite{CQ}%
\begin{eqnarray*}
&&\mathbf{W}_{H}^{3}(0,t)^{i,j,k} \\
&=&\int_{0}^{1}\left[ \int_{0}^{1}\left[ \int_{0}^{1}I_{0,t}^{K_{H}}\left(
I_{0,.}^{K_{H}}\left( K_{H}\left( .,u\right) \right) \left( v\right) \right)
\left( r\right) dB^{i}(u)\right] dB^{j}(v)\right] dB^{k}(r) \\
&&+\frac{\delta _{i,j}}{2}\int_{0}^{1}I_{0,t}^{K_{H}}\left( \mathbb{E}\left( 
\mathbf{W}_{H}^{1}(0,.)^{k}\right) ^{2}\right) \left( z\right) dB^{k}\left(
z\right) \\
&&+\frac{\delta _{k,j}}{2}\int_{0}^{1}I_{0,t}^{K_{H}}\left( \mathbb{E}\left( 
\mathbf{W}_{H}^{1}(.,t)^{i}\right) ^{2}\right) \left( u\right) dB^{i}\left(
u\right) \\
&&+\delta _{i,k}\int_{0}^{1}I_{0,t}^{K_{H}}\left( \mathbb{E}\left( \mathbf{W}%
_{H}^{1}(0,.)^{j}\mathbf{W}_{H}^{1}(.,t)^{j}\right) \right) \left( v\right)
dB^{j}\left( v\right)
\end{eqnarray*}%
The stochastic integrals with respect to the Brownian motion $B$ are
Skorokhod integrals. We then define for $s<t$%
\begin{equation*}
\mathbf{W}_{H}(s,t)=\mathbf{W}_{H}(0,s)^{-1}\otimes \mathbf{W}_{H}(0,t).
\end{equation*}
\end{definition}

\begin{remark}
\label{norm-niveau2}Using scaling property of fractional Brownian motion and
the fact that it has stationnary increment, we have for $a$ small enough%
\begin{equation*}
\sup_{s<t}{\mathbb{E}}\left( \exp a\frac{\left\Vert \mathbf{W}%
_{H}(s,t)\right\Vert ^{2}}{\left( t-s\right) ^{2H}}\right) <\infty
\end{equation*}%
Then using\textbf{\ }Garsia-Rudomich-Rumsey lemma as in \cite{FV}, it is
easy to see that%
\begin{equation*}
{\mathbb{E}}\exp \left( a\left\Vert {\mathbf{W}}_{H}\right\Vert _{\omega
,p}^{2}\right) <\infty
\end{equation*}
for some positive $a$ small enough. In particular $\left\Vert \mathbf{W}%
_{H}\right\Vert _{\omega ,p}$ belongs to $L^{q}$ for all $q\geq 1$.
\end{remark}

\subsubsection{Good Sequence for Fractional Brownian motion}

We fix a sequence $D^{n}$ of subdivision which steps size tends to $0$. We
first start by slightly extending the result of \cite{CQ}: we prove that
piecewise linear approximation of the fractional Brownian motion, lifted to
a (smooth) rough path converges to the lift of the fractional Brownian
motion constructed in \cite{CQ}. Our technique is quite different, and our
convergence stronger. Recall the notation of \ref{Dlinear}.

\begin{theorem}
\label{CQ}In the $d_{\omega ,p}$-topology, $S\left( W_{H}^{D^{n}}\right) $
converge to $\mathbf{W}_{H}$ in $L^{q}$ and in probability.\newline
\end{theorem}

\begin{proof}
For flexibel notation, we agree that%
\begin{equation*}
\mathbf{W}_{H}\left( t\right) \equiv \mathbf{W}_{t},\,\ \mathbf{W}_{H}\left(
s,t\right) \equiv \mathbf{W}_{s,t}\,.
\end{equation*}%
We fix a subdivision $D=(0\leq t_{1}<\cdots <t_{|D|}\leq 1)$ of $\left[ 0,1%
\right] .$ Similarly as for the Brownian motion, it is enough to prove that,
for $p^{\prime }>1/H,$ there exists $\mu >0$ such that for all $s$ and $t$
in $D,$ we have%
\begin{equation*}
\left\Vert \sup_{0\leq s<t\leq 1}\frac{d\left( S\left( W^{D}\right) _{s,t},%
\mathbf{W}_{s,t}\right) }{\left\vert t-s\right\vert ^{1/p}}\right\Vert
_{L^{q}}\leq C_{q,\mu }\left\vert D\right\vert ^{\mu }.
\end{equation*}%
Working just as in the case of the Brownian motion, we see that it is enough
to prove that for all $s=t_{m},t=t_{n}\in D$,%
\begin{eqnarray}
&&\left\Vert \int_{s}^{t}W_{s,u}^{D}\otimes dW_{u}^{D}-\mathbf{W}%
_{s,t}^{2}\right\Vert _{L^{2}}  \label{deg2} \\
&\leq &C_{q,\mu }\left\vert D\right\vert ^{\mu }\left\vert t-s\right\vert
^{2/p^{\prime }}, \\
&&\left\Vert \int_{s<u_{1}<u_{2}<u_{3}<t}dW_{u_{1}}^{D}\otimes
dW_{u_{2}}^{D}\otimes dW_{u_{3}}^{D}-\mathbf{W}_{s,t}^{3}\right\Vert _{L^{2}}
\label{deg3} \\
&\leq &C_{q,\mu }\left\vert D\right\vert ^{\mu }\left\vert t-s\right\vert
^{3/p^{\prime }}.
\end{eqnarray}

We will divide the proof into two parts, one corresponding to the second
level, and the other one the third level.

\textit{1st part: }A simple computation leads to%
\begin{equation*}
\int_{s}^{t}W_{s,u}^{D}\otimes dW_{u}^{D}-\mathbf{W}_{s,t}^{2}=-%
\sum_{l=m}^{n-1}\left( {\mathbf{W}}^{2}(t_{l},t_{l+1})-\frac{1}{2}%
W(t_{l},t_{l+1})^{\otimes 2}\right) .
\end{equation*}%
Therefore, to prove (\ref{deg2}), we just need to prove that%
\begin{equation*}
{\mathbb{E}}\left( \left\vert \sum_{l=m}^{n-1}\mathbf{A}_{t_{l},t_{l+1}}%
\right\vert ^{2}\right) \leq C_{q,\mu }\left\vert D\right\vert ^{\mu
}\left\vert t_{n}-t_{m}\right\vert ^{4/p\prime }
\end{equation*}%
where $\mathbf{A}_{s,t}={\mathbf{W}}^{2}(s,t)-\frac{1}{2}W(s,t)^{\otimes 2}$%
, which is the Levy area of the fractional Brownian motion. We are going to
decompose furthermore%
\begin{equation*}
\left( \mathbf{A}_{s,t}\right) ^{i,j}={\mathbf{W}}^{2}(s,t)^{i,j}-{\mathbf{W}%
}^{2}(s,t)^{j,i}
\end{equation*}%
By definition of ${\mathbf{W}}_{H}^{2}$ (or see \cite{CQ}), that for $i\neq
j $,%
\begin{align*}
{\mathbf{W}}^{2}(s,t)^{i,j}&
=\int_{0}^{1}K_{H}(t,u)W^{i}(s,u)1_{]s,t[}\left( u\right) dB_{u}^{j} \\
& +\int_{0}^{1}\left( \int_{u}^{t}W^{i}(u,r)\partial K_{H}(r,u)dr\right)
1_{]s,t[}\left( u\right) dB_{u}^{j} \\
& +\int_{0}^{1}\left( \int_{s}^{t}W^{i}(s,r)\partial K_{H}(r,u)dr\right)
1_{]0,s[}\left( u\right) dB_{u}^{j}.
\end{align*}%
Writing $\left( \mathbf{A}_{s,t}\right) ^{i,j}=\left( \mathbf{A}%
_{s,t}^{1}\right) ^{i,j}+\left( \mathbf{A}_{s,t}^{2}\right) ^{i,j}-\left( 
\mathbf{A}_{s,t}^{1}\right) ^{j,i}-\left( \mathbf{A}_{s,t}^{2}\right)
^{j,i}, $ where%
\begin{eqnarray*}
\left( \mathbf{A}_{s,t}^{1}\right) ^{i,j} &=&1_{i\neq j}\int_{0}^{1}1_{\left]
s,t\right[ }\left( u\right) \left( K_{H}(t,u)W^{i}(s,u)+\int_{u}^{t}\partial
K_{H}\left( r,u\right) W^{i}\left( u,r\right) dr\right) dB_{u}^{j} \\
\left( \mathbf{A}_{s,t}^{2}\right) ^{i,j} &=&1_{i\neq
j}\int_{0}^{1}1_{[0,s[}\left( u\right) \left( \int_{s}^{t}\partial
K_{H}(r,u)W^{i}(s,r)dr\right) dB_{u}^{j}.
\end{eqnarray*}%
First observe that%
\begin{eqnarray*}
\left\vert \sum_{l=m}^{n-1}\mathbf{A}_{t_{l},t_{l+1}}\right\vert ^{2}
&=&\sum_{i,j=1}^{d}\left\vert \sum_{l=m}^{n-1}\left( \mathbf{A}%
_{t_{l},t_{l+1}}\right) ^{i,j}\right\vert ^{2} \\
&\leq &C\sum_{i,j}\left\vert \sum_{l=m}^{n-1}\left( \mathbf{A}%
_{t_{l},t_{l+1}}^{1}\right) ^{i,j}\right\vert ^{2}+C\sum_{i,j}\left\vert
\sum_{l=m}^{n-1}\left( \mathbf{A}_{t_{l},t_{l+1}}^{2}\right)
^{i,j}\right\vert ^{2}.
\end{eqnarray*}

The lemmae in the next appendix (lemmae \ref{1lem},\ref{2lem}) show that for
all $s^{\prime }<t^{\prime }\leq s<t$%
\begin{equation*}
\left\vert {\mathbb{E}}\left( \left( \mathbf{A}_{s^{\prime },t^{\prime
}}^{1}\right) ^{i,j}\left( \mathbf{A}_{s,t}^{1}\right) ^{i,j}\right)
\right\vert +\left\vert {\mathbb{E}}\left( \left( \mathbf{A}_{s^{\prime
},t^{\prime }}^{2}\right) ^{i,j}\left( \mathbf{A}_{s,t}^{2}\right)
^{i,j}\right) \right\vert \leq {\mathbb{E}}(W^{i}(s,t)W^{i}(s^{\prime
},t^{\prime }))^{2}.
\end{equation*}%
From the scaling property of the Fractional Brownian Motion, we have%
\begin{equation*}
{\mathbb{E}}\left( \left\vert \mathbf{A}_{s,t}\right\vert ^{2}\right) \leq
C\left\vert t-s\right\vert ^{4H}=C{\mathbb{E}}\left( \left\vert
W^{i}(s,t)\right\vert ^{2}\right) ^{2}
\end{equation*}%
With lemma \ref{2lemcopy}, we therefore obtain for $s<t,$%
\begin{equation*}
{\mathbb{E}}\left( \left\vert \left( \mathbf{A}_{s,t}^{1}\right)
^{i,j}\right\vert ^{2}\right) +{\mathbb{E}}\left( \left\vert \left( \mathbf{A%
}_{s,t}^{2}\right) ^{i,j}\right\vert ^{2}\right) \leq C{\mathbb{E}}\left(
\left\vert W_{H}^{i}(s,t)\right\vert ^{2}\right) ^{2}.
\end{equation*}%
Hence, we obtain that%
\begin{equation*}
{\mathbb{E}}\left( \left\vert \sum_{l=m}^{n-1}\mathbf{A}_{t_{l},t_{l+1}}%
\right\vert ^{2}\right) \leq C\sum_{k,l=m}^{n-1}{\mathbb{E}}\left(
W^{1}(t_{k},t_{k+1})W^{1}(t_{l},t_{l+1})\right) ^{2}.
\end{equation*}%
Therefore, by lemma \ref{finallemma}, we obtain%
\begin{equation*}
{\mathbb{E}}\left( \left\vert \sum_{l=m}^{n-1}\mathbf{A}_{t_{l},t_{l+1}}%
\right\vert ^{2}\right) \leq C\left\vert D\right\vert ^{4H-4/p^{\prime
}}\left\vert t_{n}-t_{m}\right\vert ^{4/p\prime },
\end{equation*}%
i.e. we have proved equation (\ref{deg2}).

\textit{2nd part:} First observe that 
\begin{equation}
\left( \mathbf{W}_{s,t}^{3}\right) ^{i,k,j}=\left( \mathbf{W}%
_{s,t}^{2}\right) ^{i,k}\left( \mathbf{W}_{s,t}^{1}\right) ^{j}-\left( 
\mathbf{W}_{s,t}^{3,i,j,k}+\mathbf{W}_{s,t}^{3,j,i,k}\right) .
\label{niveau2-1}
\end{equation}%
Therefore, using twice (\ref{niveau2-1}), for $(i,j,j)$ and $%
(i,j,k)=(j,i,j), $\newline
$\left( \mathbf{W}_{s,t}^{3}\right) ^{i,j,i}=\left( \mathbf{W}%
_{s,t}^{2}\right) ^{i,j}\left( \mathbf{W}_{s,t}^{1}\right) ^{i}-2\mathbf{W}%
_{s,t}^{i,i,j}$ and%
\begin{eqnarray*}
2\left( \mathbf{W}_{s,t}^{3}\right) ^{i,j,j} &=&\left( \mathbf{W}%
_{s,t}^{2}\right) ^{i,j}\left( \mathbf{W}_{s,t}^{1}\right) ^{j}-\mathbf{W}%
_{s,t}^{j,i,j} \\
&=&\left\{ \left( \mathbf{W}_{s,t}^{2}\right) ^{i,j}-\left( \mathbf{W}%
_{s,t}^{2}\right) ^{j,i}\right\} \left( \mathbf{W}_{s,t}^{1}\right) ^{j}+2%
\mathbf{W}_{s,t}^{j,j,i}.
\end{eqnarray*}%
In particular, as, for $s<t\in D,$ 
\begin{eqnarray*}
\left\Vert \left\{ \mathbf{W}_{s,t}^{2}-S\left( W^{D}\right)
_{s,t}^{2}\right\} \otimes \mathbf{W}_{s,t}^{1}\right\Vert _{L^{2}} &\leq
&\left\Vert \left\{ \mathbf{W}_{s,t}^{2}-S\left( W^{D}\right)
_{s,t}^{2}\right\} \right\Vert _{L^{4}}\left\Vert \mathbf{W}%
_{s,t}^{1}\right\Vert _{L^{4}} \\
&\leq &C\left\vert D\right\vert ^{4H-4/p^{\prime }}\left\vert t-s\right\vert
^{2/p\prime +H},
\end{eqnarray*}%
we see that we just need to prove that for all $i,j,k$ with $k\notin \left\{
i,j\right\} ,$%
\begin{equation}
\left\Vert \mathbf{W}_{s,t}^{3,i,j,k}+\mathbf{W}_{s,t}^{3,j,i,k}-\left(
S\left( W^{D}\right) _{s,t}^{3}\right) ^{i,j,k}-\left( S\left( W^{D}\right)
_{s,t}^{3}\right) ^{j,i,k}\right\Vert _{L^{2}}\leq C_{q,\mu }\left\vert
D\right\vert ^{\mu }\left\vert t-s\right\vert ^{3/p^{\prime }}.
\label{degree3}
\end{equation}%
The multiplicative (or Chen) property of iterated integrals quickly gives
that $\mathbf{W}_{t_{m},t_{n}}^{3,i,j,k}+\mathbf{W}_{t_{m},t_{n}}^{3,j,i,k}-%
\left( S\left( W^{D}\right) _{t_{m},t_{n}}^{3}\right) ^{i,j,k}-\left(
S\left( W^{D}\right) _{t_{m},t_{n}}^{3}\right) ^{j,i,k}$ is equal to%
\begin{eqnarray*}
&&\sum_{l=m}^{n-1}\mathbf{W}_{t_{l},t_{l+1}}^{3,i,j,k}+\mathbf{W}%
_{t_{l},t_{l+1}}^{3,j,i,k}-\left( S\left( W^{D}\right)
_{t_{l},t_{l+1}}^{3}\right) ^{i,j,k}-\left( S\left( W^{D}\right)
_{t_{l},t_{l+1}}^{3}\right) ^{j,i,k} \\
&&-\sum_{l=m}^{n-1}W_{t_{m},t_{l}}^{i}\left( \frac{1}{2}%
W_{t_{l},t_{l+1}}^{j}W_{t_{l},t_{l+1}}^{k}-\mathbf{W}%
_{t_{l},t_{l+1}}^{2,j,k}\right) \\
&&-\sum_{l=m}^{n-1}W_{t_{m},t_{l}}^{j}\left( \frac{1}{2}%
W_{t_{l},t_{l+1}}^{i}W_{t_{l},t_{l+1}}^{k}-\mathbf{W}%
_{t_{l},t_{l+1}}^{2,i,k}\right) .
\end{eqnarray*}%
We first start by bounding the $L^{2}$ norm of the second line (we bound the
third line in just the same way):%
\begin{multline*}
\left\Vert \sum_{l=m}^{n-1}W_{t_{m},t_{l}}^{i}\left( \frac{1}{2}%
W_{t_{l},t_{l+1}}^{j}W_{t_{l},t_{l+1}}^{k}-\mathbf{W}%
_{t_{l},t_{l+1}}^{2,j,k}\right) \right\Vert _{L^{2}}=\frac{1}{2}\left\Vert
\sum_{l=m}^{n-1}W_{t_{m},t_{l}}^{i}W_{t_{l},t_{l+1}}^{,j}W_{t_{l},t_{l+1}}^{k}\right\Vert _{L^{2}}
\\
+\sum_{\alpha =1}^{2}\left\Vert \sum_{l=m}^{n-1}W_{t_{m},t_{l}}^{i}\left( 
\mathbf{A}_{t_{l},t_{l+1}}^{\alpha ,1}\right) ^{j,k}\right\Vert _{L^{2}},
\end{multline*}%
where $A^{\alpha },~~\alpha =1,2$ are defined earlier in this proof. Because 
$k\neq j,$ using the Gaussian formula (\ref{Gaussian}) 
\begin{eqnarray*}
&&{\mathbb{E}}\left( \left\vert
\sum_{l=m}^{n-1}W_{t_{m},t_{l}}^{i}W_{t_{l},t_{l+1}}^{,j}W_{t_{l},t_{l+1}}^{k}\right\vert ^{2}\right)
\\
&=&\sum_{l,l^{\prime }=m}^{n-1}{\mathbb{E}}\left(
W_{t_{m},t_{l}}^{i}W_{t_{l},t_{l+1}}^{j}W_{t_{m},t_{l^{\prime
}}}^{i}W_{t_{l^{\prime }},t_{l^{\prime }+1}}^{,j}\right) {\mathbb{E}}\left(
W_{t_{l},t_{l+1}}^{k}W_{t_{l^{\prime }},t_{l^{\prime }+1}}^{k}\right) \\
&\leq &\sum_{l,l^{\prime }=m}^{n-1}\left\Vert W_{t_{m},t_{l}}^{i}\right\Vert
_{L^{4}}\left\Vert W_{t_{l},t_{l+1}}^{j}\right\Vert _{L^{4}}\left\Vert
W_{t_{m},t_{l^{\prime }}}^{i}\right\Vert _{L^{4}}\left\Vert W_{t_{l^{\prime
}},t_{l^{\prime }+1}}^{,j}\right\Vert _{L^{4}}{\mathbb{E}}\left(
W_{t_{l},t_{l+1}}^{k}W_{t_{l^{\prime }},t_{l^{\prime }+1}}^{k}\right) \\
&\leq &C\left\vert t_{n}-t_{m}\right\vert ^{2H}\sum_{l,l^{\prime
}=m}^{n-1}\left\Vert W_{t_{l},t_{l+1}}^{j}\right\Vert _{L^{2}}\left\Vert
W_{t_{l^{\prime }},t_{l^{\prime }+1}}^{j}\right\Vert _{L^{2}}{\mathbb{E}}%
\left( W_{t_{l},t_{l+1}}^{k}W_{t_{l^{\prime }},t_{l^{\prime }+1}}^{k}\right)
\\
&\leq &C\left\vert t_{n}-t_{m}\right\vert ^{6/p^{\prime }}\left\vert
D\right\vert ^{\mu }\text{ \ by lemma \ref{finallemma}.}
\end{eqnarray*}%
In lemma \ref{2lemdegree3}, we proved that%
\begin{align*}
& \left\Vert \sum_{l=m}^{n-1}W_{t_{m},t_{l}}^{i}\left( \mathbf{A}%
_{t_{l},t_{l+1}}^{2}\right) ^{j,k}\right\Vert _{L^{2}}^{2} \\
& \leq \sum_{l,l^{\prime }=m}^{n-1}{\mathbb{E}}\left(
W_{t_{m},t_{l}}^{i}W_{t_{l},t_{l+1}}^{j}W_{t_{m},t_{l^{\prime
}}}^{i}W_{t_{l^{\prime }},t_{l^{\prime }+1}}^{j}\right) {\mathbb{E}}\left(
W_{t_{l},t_{l+1}}^{k}W_{t_{l^{\prime }},t_{l^{\prime }+1}}^{k}\right) \\
& \leq \left\vert t_{n}-t_{m}\right\vert ^{6/p^{\prime }}\left\vert
D\right\vert ^{\mu }\text{ \ as above.}
\end{align*}%
Using the fact that $k\notin \left\{ i,j\right\} $ and lemma \ref{1lem},%
\begin{eqnarray*}
&&\left\Vert \sum_{l=m}^{n-1}W_{t_{m},t_{l}}^{i}\left( \mathbf{A}%
_{t_{l},t_{l+1}}^{1}\right) ^{j,k}\right\Vert _{L^{2}} \\
&=&\left( \sum_{l=m}^{n-1}{\mathbb{E}}\left( \left\vert
W_{t_{m},t_{l}}^{i}\left( \mathbf{A}_{t_{l},t_{l+1}}^{1}\right)
^{j,k}\right\vert ^{2}\right) \right) ^{1/2} \\
&\leq &\left( \sum_{l=m}^{n-1}{\mathbb{E}}\left( \left\vert \left( \mathbf{A}%
_{t_{l},t_{l+1}}^{1}\right) ^{j,k}\right\vert ^{4}\right) ^{1/2}{\mathbb{E}}%
\left( \left\vert W_{t_{m},t_{l}}^{i}\right\vert ^{4}\right) ^{1/2}\right)
^{1/2} \\
&\leq &\left( \left\vert t_{n}-t_{m}\right\vert
^{2H}\sum_{l=m}^{n-1}\left\vert t_{l+1}-t_{l}\right\vert ^{4H}\right) ^{1/2}
\\
&\leq &\left( C\left\vert t_{n}-t_{m}\right\vert ^{2H}\left\vert
t_{n}-t_{m}\right\vert ^{4/p^{\prime }}|D|^{4h-4/p^{\prime }}\right) ^{1/2}
\\
&\leq &C_{q,\mu }\left\vert D\right\vert ^{\mu }\left\vert t-s\right\vert
^{3/p^{\prime }}.
\end{eqnarray*}%
Therefore, to prove equation (\ref{degree3}), it remains to prove that for $%
k\notin \left\{ i,j\right\} $%
\begin{equation*}
\left\Vert \sum_{l=m}^{n-1}\mathbf{W}_{t_{l},t_{l+1}}^{3,i,j,k}-\left(
S\left( W^{D}\right) _{t_{l},t_{l+1}}^{3}\right) ^{i,j,k}\right\Vert
_{L^{2}}\leq C_{q,\mu }\left\vert D\right\vert ^{\mu }\left\vert
t-s\right\vert ^{3/p^{\prime }}.
\end{equation*}%
First observe that $\left( S\left( W^{D}\right) _{t_{l},t_{l+1}}^{3}\right)
^{i,j,k}=\frac{1}{6}%
W_{t_{l},t_{l+1}}^{i}W_{t_{l},t_{l+1}}^{j}W_{t_{l},t_{l+1}}^{k}$ and that (
using the Gaussian formula (\ref{Gaussian}))%
\begin{eqnarray*}
&&\left\Vert
\sum_{l=m}^{n-1}W_{t_{l},t_{l+1}}^{i}W_{t_{l},t_{l+1}}^{j}W_{t_{l},t_{l+1}}^{k}\right\Vert _{L^{2}}^{2}
\\
&\leq &\sum_{l,l^{\prime }=m}^{n-1}{\mathbb{E}}\left(
W_{t_{l},t_{l+1}}^{i}W_{t_{l},t_{l+1}}^{j}W_{t_{l^{\prime }},t_{l^{\prime
}+1}}^{i}W_{t_{l^{\prime }},t_{l^{\prime }+1}}^{j}\right) {\mathbb{E}}\left(
W_{t_{l},t_{l+1}}^{k}W_{t_{l^{\prime }},t_{l^{\prime }+1}}^{k}\right) \\
&\leq &C\left\vert t_{n}-t_{m}\right\vert ^{6/p^{\prime }}\left\vert
D\right\vert ^{\mu }\text{ \ by lemma \ref{finallemma}.}
\end{eqnarray*}%
Hence, we need to prove that%
\begin{equation*}
\left\Vert \sum_{l=m}^{n-1}\mathbf{W}_{t_{l},t_{l+1}}^{3,i,j,k}\right\Vert
_{L^{2}}\leq C_{q,\mu }\left\vert D\right\vert ^{\mu }\left\vert
t-s\right\vert ^{3/p^{\prime }}.
\end{equation*}%
Recall that for $i,j,k$ as above, we have%
\begin{align*}
{\mathbf{W}}^{3}(s,t)^{i,j,k}& =\int_{0}^{1}K_{H}(t,u)\mathbf{W}%
^{2}(s,u)^{i,j}1_{]s,t[}\left( u\right) dB_{u}^{k} \\
& +\int_{0}^{1}\left( \int_{u}^{t}\mathbf{W}^{2}(u,r)^{i,j}\partial
K_{H}(r,u)dr\right) 1_{]s,t[}\left( u\right) dB_{u}^{k} \\
& +\int_{0}^{1}\left( \mathbf{W}^{1}(s,u)^{i}\int_{u}^{t}\mathbf{W}%
^{1}(u,r)^{j}\partial K_{H}(r,u)dr\right) 1_{]s,t[}\left( u\right) dB_{u}^{k}
\\
& +\int_{0}^{1}\left( \int_{s}^{t}\mathbf{W}^{2}(s,r)^{i,j}\partial
K_{H}(r,u)dr\right) 1_{]0,s[}\left( u\right) dB_{u}^{k}.
\end{align*}

Since $\mathbf{W}^{2}(s,t)=\mathbf{A}_{s,t}^{1}+\mathbf{A}_{s,t}^{2},$ we let%
\begin{align*}
\mathbf{G}_{H}^{1}(s,t)^{i,j,k}& =\int_{0}^{1}K_{H}(t,u)\mathbf{W}%
^{2}(s,u)^{i,j}1_{]s,t[}\left( u\right) dB_{u}^{k} \\
& +\int_{0}^{1}\left( \int_{u}^{t}\mathbf{W}^{2}(u,r)^{i,j}\partial
K_{H}(r,u)dr\right) 1_{]s,t[}\left( u\right) dB_{u}^{k} \\
& +\int_{0}^{1}\left( \mathbf{W}^{1}(s,u)^{i}\int_{u}^{t}\mathbf{W}%
^{1}(u,r)^{j}\partial K_{H}(r,u)dr\right) 1_{]s,t[}\left( u\right) dB_{u}^{k}
\end{align*}%
and for $\alpha \in \left\{ 1,2\right\} $%
\begin{equation*}
\mathbf{G}_{H}^{2,\alpha }(s,t)^{i,j,k}=\int_{0}^{1}\left(
\int_{s}^{t}\left( \mathbf{A}_{s,r}^{\alpha }\right) ^{i,j}\partial
K_{H}(r,u)dr\right) 1_{]0,s[}\left( u\right) dB_{u}^{k},
\end{equation*}%
so that ${\mathbf{W}}_{H}^{3}(s,t)^{i,j,k}=\mathbf{G}_{H}^{1}(s,t)^{i,j,k}+%
\mathbf{G}_{H}^{2,1}(s,t)^{i,j,k}+\mathbf{G}_{H}^{2,2}(s,t)^{i,j,k}.$
Arguments of lemmae \ref{1lem} and \ref{2lem} leads to, for $\alpha \in
\left\{ 1,2\right\} $ and $s^{\prime }<t^{\prime }\leq s<t$ or $\left(
s^{\prime }<t^{\prime }\right) =\left( s<t\right) ,$%
\begin{equation*}
{\mathbb{E}}\left( \mathbf{G}_{H}^{2,\alpha }(s,t)^{i,j,k}\mathbf{G}%
_{H}^{2,\alpha }(s^{\prime },t^{\prime })^{i,j,k}\right) \leq -C{\mathbb{E}}%
(W_{H}^{i}(s,t)W_{H}^{i}(s^{\prime },t^{\prime }))^{3}.
\end{equation*}%
Therefore, working as above, we get%
\begin{equation*}
\left\Vert \sum_{l=m}^{n-1}\mathbf{G}_{H}^{2,\alpha
}(t_{l},t_{l+1})^{i,j,k}\right\Vert _{L^{2}}\leq C_{q,\mu }\left\vert
D\right\vert ^{\mu }\left\vert t-s\right\vert ^{3/p^{\prime }}.
\end{equation*}%
Then, lemma \ref{1lem} gives that%
\begin{equation*}
{\mathbb{E}}\left( \mathbf{G}_{H}^{1}(s,t)^{i,j,k}\mathbf{G}%
_{H}^{1}(s^{\prime },t^{\prime })^{i,j,k}\right) =0
\end{equation*}%
for $s^{\prime }<t^{\prime }\leq s<t$ and scaling property of the fractional
Brownian motion gives%
\begin{equation*}
\left\Vert \mathbf{G}_{H}^{1}(s,t)^{i,j,k}\right\Vert _{L^{2}}=C\left\vert
t-s\right\vert ^{3H},
\end{equation*}%
so that we once again easily find%
\begin{equation*}
\left\Vert \sum_{l=m}^{n-1}\mathbf{G}_{H}^{2,\alpha
}(t_{l},t_{l+1})^{i,j,k}\right\Vert _{L^{2}}\leq C_{q,\mu }\left\vert
D\right\vert ^{\mu }\left\vert t-s\right\vert ^{3/p^{\prime }}.
\end{equation*}%
That concludes the proof.
\end{proof}

We go one step further than in the paper \cite{CQ}, by proving the following
theorem:

\begin{theorem}
\label{Bngoodmbf}When $n$ tends to $\infty ,$ $d_{\omega ,p}\left( S^{\prime
}\left( W_{H}^{D^{n}},\mathbf{W}\right) ,S^{\prime \prime }\left( \mathbf{W},%
\mathbf{W}\right) \right) $ converges to $0$ in $L^{q}$, $q\geq 1$ and in
probability.\newline
If $D^{n}=\left( \frac{k}{2^{n}},0\leq k\leq 2^{n}\right) $, the convergence
also holds almost surely, i.e. $W_{H}^{n}=W_{H}^{D^{n}}$ is almost surely a
good $p$-rough path sequence associated to $\mathbf{W}_{H}$.
\end{theorem}

We first start by two lemmae.

\begin{lemma}
\label{leadtothemainth}Let $D$ be a subdivision of $[0,1]$. Then, there
exists $\mu >0$ and $C_{\mu }<\infty $ such that for all $s,t\in D,$%
\begin{equation}
\left\Vert \int_{s}^{t}W\left( s,u\right) \otimes dW^{D}\left( u\right) -%
\mathbf{W}_{s,t}^{2}\right\Vert _{L^{2}}\leq C_{\mu }\left\vert D\right\vert
^{\mu }\left\vert t-s\right\vert ^{2/p^{\prime }}.  \label{sansDFBM}
\end{equation}
\end{lemma}

\begin{proof}
The subdivision $D=(0\leq t_{1}<\cdots <t_{|D|}\leq 1)$ is fixed. As%
\begin{equation*}
\left\Vert \int_{s}^{t}W^{D}\left( s,u\right) \otimes dW^{D}\left( u\right) -%
\mathbf{W}_{s,t}^{2}\right\Vert _{L^{2}}\leq C_{\mu }\left\vert D\right\vert
^{\mu }\left\vert t-s\right\vert ^{4/p^{\prime }},
\end{equation*}%
we just need to prove that for all $0\leq m<n\leq \left\vert D\right\vert $,%
\begin{equation*}
{\mathbb{E}}\left( \left\vert \int_{t_{m}}^{t_{n}}\left(
W^{D}(t_{m},u)-W\left( t_{m},u\right) \right) \otimes dW^{D}\left( u\right)
\right\vert ^{2}\right) \leq C\left\vert D\right\vert ^{\mu }\left\vert
t_{n}-t_{m}\right\vert ^{4/p\prime },
\end{equation*}%
i.e. that%
\begin{equation*}
{\mathbb{E}}\left( \left\vert \sum_{j=m}^{n-1}\int_{t_{j}}^{t_{j+1}}\left(
W^{D}(t_{j},u)-W\left( t_{j},u\right) \right) \otimes dW^{D}\left( u\right)
\right\vert ^{2}\right) \leq C\left\vert D\right\vert ^{\mu }\left\vert
t_{n}-t_{m}\right\vert ^{4/p\prime }
\end{equation*}%
For $s<t$ , we let 
\begin{equation*}
\mathbf{Y}_{s,t}=\frac{1}{|t-s|^{2}}\int_{s}^{t}[(u-s)W(u,t)-(t-u)W(s,u)]du%
\otimes {\mathbf{W}}^{1}(s,t).
\end{equation*}%
As for $s,t$ are two consecutive points in $D,$ $\mathbf{Y}%
_{s,t}=\int_{s}^{t}\left( W^{D}(s,u)-W\left( s,u\right) \right) \otimes
dW^{D}\left( u\right) $, to prove (\ref{sansDFBM}), we just need to prove
that%
\begin{equation}
{\mathbb{E}}\left( \left\vert \sum_{l=m}^{n-1}\mathbf{Y}_{t_{l},t_{l+1}}%
\right\vert ^{2}\right) \leq C\left\vert D\right\vert ^{\theta }\left\vert
t_{n}-t_{m}\right\vert ^{2/p\prime }.  \label{eq2}
\end{equation}

Lemma \ref{3lem} that for all $s^{\prime }<t^{\prime }\leq s<t$%
\begin{equation*}
{\mathbb{E}}\left( \left( \mathbf{Y}_{s^{\prime },t^{\prime }}\right)
^{i,j}\left( \mathbf{Y}_{s,t}\right) ^{i,j}\right) \leq C{\mathbb{E}}%
(W^{i}(s,t)W^{i}(s^{\prime },t^{\prime }))^{2}.
\end{equation*}%
From the scaling property of the fractional Brownian Motion, we have%
\begin{equation*}
{\mathbb{E}}\left( \left\vert \mathbf{Y}_{s,t}\right\vert ^{2}\right) \leq
C\left\vert t-s\right\vert ^{4H}=C{\mathbb{E}}\left( \left\vert
W^{i}(s,t)\right\vert ^{2}\right) ^{2}
\end{equation*}%
Hence, we obtain%
\begin{equation*}
{\mathbb{E}}\left( \left\vert \sum_{l=m}^{n-1}\mathbf{Y}_{t_{l},t_{l+1}}%
\right\vert ^{2}\right) \leq C\sum_{k,l=m}^{n-1}E\left(
W^{1}(t_{k},t_{k+1})W^{1}(t_{l},t_{l+1})\right) ^{2}.
\end{equation*}%
which gives by lemma \ref{finallemma}%
\begin{equation*}
{\mathbb{E}}\left( \left\vert \sum_{l=m}^{n-1}\mathbf{Y}_{t_{l},t_{l+1}}%
\right\vert ^{2}\right) \leq C\left\vert D\right\vert ^{4H-4/p^{\prime
}}\left\vert t_{n}-t_{m}\right\vert ^{4/p\prime }.
\end{equation*}
\end{proof}

\begin{lemma}
\label{leadtothemainthdegree3}Let $D$ be a subdivision of $[0,1]$. Then, for
all $q\geq 1,$ there exists $\mu >0$ and $C_{q,\mu }<\infty $ such that for
all $s<t\in D,$

\begin{align*}
& \left\Vert \int_{s<u_{1}<u_{2}<u_{3}<t}\left\{ dW\left( u_{1}\right)
-dW^{D}\left( u_{1}\right) \right\} \otimes dW^{D}\left( u_{2}\right)
\otimes dW^{D}\left( u_{3}\right) \right\Vert _{L^{q}} \\
& \leq C_{q,\mu }\left\vert D\right\vert ^{\mu }\left\vert t-s\right\vert
^{3/p^{\prime }},
\end{align*}%
and%
\begin{equation*}
\left\Vert \int_{s<u<t}\left\{ \mathbf{W}^{2}\left( s,u\right) -S\left(
W^{D}\right) _{s,u}^{2}\right\} \otimes dW^{D}\left( u\right) \right\Vert
_{L^{q}}\leq C_{q,\mu }\left\vert D\right\vert ^{\mu }\left\vert
t-s\right\vert ^{3/p^{\prime }}
\end{equation*}
\end{lemma}

\begin{proof}
The subdivision $D=(0\leq t_{1}<\cdots <t_{|D|}\leq 1)$ is fixed. With
similar techniques as in theorem \ref{CQ}, we see that to obtain the first
inequality, we just need to prove that%
\begin{gather}
{\mathbb{E}}\left( \left\vert
\sum_{l=m}^{n-1}\int_{t_{l}<u_{1}<u_{2}<u_{3}<t_{l+1}}\left\{ dW\left(
u_{1}\right) -dW^{D}\left( u_{1}\right) \right\} \otimes dW^{D}\left(
u_{2}\right) \otimes dW^{D}\left( u_{3}\right) \right\vert ^{2}\right) 
\notag \\
\leq C\left\vert D\right\vert ^{\mu }\left\vert t_{n}-t_{m}\right\vert
^{6/p\prime },  \label{2D'}
\end{gather}

while to obtain the second one, it is enough to prove that%
\begin{equation}
{\mathbb{E}}\left( \left\vert \sum_{l=m}^{n-1}\int_{t_{l}}^{t_{l+1}}\mathbf{W%
}_{t_{l},u}^{2,i,j}dW^{D,k}\left( u\right) \right\vert ^{2}\right) \leq
C\left\vert D\right\vert ^{\mu }\left\vert t_{n}-t_{m}\right\vert
^{6/p\prime }.  \label{1D'}
\end{equation}

We first start by proving equation (\ref{2D'}). We define%
\begin{equation*}
\mathbf{\tilde{Y}}_{s,t}=\frac{1}{|t-s|^{3}}%
\int_{s<u_{1}<u_{2}<t}[(u_{1}-s)W(u_{1},t)-(t-u_{1})W(s,u_{1})]du_{1}du_{2}%
\otimes \left( {\mathbf{W}}^{1}(s,t)\right) ^{\otimes 2}
\end{equation*}%
so that $\mathbf{\tilde{Y}}_{s,t}=\int_{s<u_{1}<u_{2}<u_{3}<t}\left\{
dW\left( u_{1}\right) -dW^{D}\left( u_{1}\right) \right\} \otimes
dW^{D}\left( u_{2}\right) \otimes dW^{D}\left( u_{3}\right) .$ From lemma %
\ref{3lembis}, we have that%
\begin{eqnarray*}
&&{\mathbb{E}}\left( \left\vert \sum_{l=m}^{n-1}\mathbf{\tilde{Y}}%
_{t_{l},t_{l+1}}\right\vert ^{2}\right) \\
&\leq &\sum_{k,l}{\mathbb{E}}\left( \left\vert
W_{t_{k},t_{k+1}}^{1}\right\vert ^{2}\right) {\mathbb{E}}\left( \left\vert
W_{t_{l},t_{l+1}}^{1}\right\vert ^{2}\right) {\mathbb{E}}\left(
W_{t_{k},t_{k+1}}^{1}W_{t_{l},t_{l+1}}^{1}\right) \\
&\leq &\left\vert t_{n}-t_{m}\right\vert ^{2h^{\prime }}\sum_{k,l}{\mathbb{E}%
}\left( \left\vert W_{t_{k},t_{k+1}}^{1}\right\vert ^{2}\right) ^{1/2}{%
\mathbb{E}}\left( \left\vert W_{t_{l},t_{l+1}}^{1}\right\vert ^{2}\right)
^{1/2}{\mathbb{E}}\left( W_{t_{k},t_{k+1}}^{1}W_{t_{l},t_{l+1}}^{1}\right) \\
&\leq &C\left\vert t_{n}-t_{m}\right\vert ^{6/p^{\prime }}\left\vert
D\right\vert ^{\mu }\text{ \ by lemma \ref{finallemma}.}
\end{eqnarray*}%
We can prove (\ref{1D'}) when $i=j,$ with a simple use of the Gaussian
formula (\ref{Gaussian}), in a similar way as the proof or equation (\ref%
{2D'}). The case $k\notin \left\{ i,j\right\} $ is not too difficult, using
the independence between $W_{H}^{D,k}$ and $\mathbf{W}_{t_{l},u}^{2,i,j}$.
So we take $k\in \left\{ i,j\right\} $, with $i\neq j$. An integration by
part shows that the cases $k=i$ and $k=j$ are equivalent. Hence, we just
need to prove equation (\ref{1D'}) for $k=i$ different than $j$. We let for $%
\alpha \in \left\{ 1,2\right\} $,%
\begin{equation*}
\mathbf{Z}_{s,t}^{\alpha }=\left( \frac{1}{t-s}\int_{s}^{t}\left( \mathbf{A}%
_{s,u}^{\alpha }\right) ^{i,j}du\right) W^{i}\left( s,t\right)
\end{equation*}%
so that for $s$ and $t$ in $D,$ we have $\int_{t_{l}<u<t_{l+1}}\mathbf{W}%
_{t_{l},u}^{2,i,j}dW_{H}^{D,k}\left( u\right) =\mathbf{Z}_{s,t}^{1}+\mathbf{Z%
}_{s,t}^{2}.$ Lemma \ref{1lem} gives ${\mathbb{E}}\left( \mathbf{Z}_{s,t}^{1}%
\mathbf{Z}_{s^{\prime },t^{\prime }}^{1}\right) =0$ whenever $s<t\leq
s^{\prime }<t^{\prime },$ while lemma \ref{forthirdlevel} gives ${\mathbb{E}}%
\left( \mathbf{Z}_{s,t}^{2}\mathbf{Z}_{s^{\prime },t^{\prime }}^{2}\right)
\leq C{\mathbb{E}}\left( W^{j}(s,t)W^{j}(s^{\prime },t^{\prime })\right) {%
\mathbb{E}}\left( \left\vert W^{j}(s,t)\right\vert ^{2}\right) {\mathbb{E}}%
\left( \left\vert W^{j}(s^{\prime },t^{\prime })\right\vert ^{2}\right) $
for all $s<t\leq s^{\prime }<t^{\prime }$ or $\left( s<t\right) =\left(
s^{\prime }<t^{\prime }\right) .$ Then, the scaling property of the
fractional Brownian motion gives ${\mathbb{E}}\left( \left\vert \mathbf{Z}%
_{s,t}^{1}\right\vert ^{2}\right) =C{\mathbb{E}}\left( \left\vert
W^{j}(s,t)\right\vert ^{2}\right) ^{3},$ so that%
\begin{eqnarray*}
&&{\mathbb{E}}\left( \left\vert \sum_{l=m}^{n-1}\int_{t_{l}<u<t_{l+1}}%
\mathbf{W}_{t_{l},u}^{2,i,j}dW^{D,k}\left( u\right) \right\vert ^{2}\right)
\\
&\leq &C\sum_{k,l}{\mathbb{E}}\left( \left\vert
W_{t_{k},t_{k+1}}^{1}\right\vert ^{2}\right) {\mathbb{E}}\left( \left\vert
W_{t_{l},t_{l+1}}^{1}\right\vert ^{2}\right) {\mathbb{E}}\left(
W_{t_{k},t_{k+1}}^{1}W_{t_{l},t_{l+1}}^{1}\right) \\
&\leq &C\left\vert t_{n}-t_{m}\right\vert ^{6/p^{\prime }}\left\vert
D\right\vert ^{\mu }\text{ \ by lemma \ref{finallemma}.}
\end{eqnarray*}
\end{proof}

These previous lemmae, a few more integrations by part, and techniques as in
theorem \ref{Bngood} easily lead to theorem \ref{Bngoodmbf}.

\subsection{$W^{1,\infty }$ paths are in the Cameron-Martin Space\label%
{plaacm}}

We check that the piece-wise linear combinations are in the Cameron-Martin
space associated to Fractional Brownian motion, $\mathcal{H}_{H}$. Since the
case of Brownian motion, that is $H=1/2$, is trivial we focus on $H<1/2.$%
Recall \cite{Due} that the integral kernel $K_{H}\left( t,s\right) $ gives
rise to a well-defined integral transform for $L^{2}\left( \left[ 0,1\right]
\right) $-functions and%
\begin{equation*}
\mathcal{H}_{H}=K_{H}\left[ L^{2}\left( \left[ 0,1\right] \right) \right]
\end{equation*}

\begin{lemma}
Let $H\in (0,1/2)$. Every function $f$ with $f\left( 0\right) =0$ and $%
\left\vert f^{\prime }\right\vert \in L^{\infty }$ is in the Cameron-Martin
space $\mathcal{H}_{H}.$
\end{lemma}

\begin{proof}
Theorem 2.1 of \cite{Due} states%
\begin{equation*}
\mathcal{H}_{H}=I_{0+}^{H+1/2}\left[ L^{2}\left( \left[ 0,1\right] \right) %
\right]
\end{equation*}%
where fractional integral operator%
\begin{equation*}
I_{0+}^{\alpha }\left[ f\right] \left( x\right) =\frac{1}{\Gamma \left(
\alpha \right) }\int_{0}^{x}f\left( t\right) \left( x-t\right) ^{\alpha
-1}dt,\text{ \ \ \ \ }x\in \left[ 0,1\right] .
\end{equation*}%
Restricting to $\alpha \in (0,1),$ its inverse is given by the fractional
derivative%
\begin{equation*}
D_{0+}^{\alpha }\left[ f\right] \left( t\right) \equiv \frac{d}{dt}%
I_{0+}^{1-\alpha }[f]\left( t\right) .
\end{equation*}%
So it will be enough to prove that for $f$ smooth and $q$ large enough,%
\begin{equation*}
\left\Vert D_{0+}^{h+1/2}\left[ f\right] \right\Vert _{L^{2}}\leq
C_{H,q}\left\Vert f^{\prime }\right\Vert _{L^{q}}.
\end{equation*}%
This is easy as 
\begin{eqnarray*}
\left\vert D_{0+}^{H+1/2}\left[ f\right] \left( t\right) \right\vert
&=&\left\vert \frac{1}{\Gamma \left( \alpha \right) }\frac{d}{dt}%
\int_{0}^{t}f\left( r\right) \left( t-r\right) ^{-H-1/2}dr\right\vert \\
&=&\frac{H+1/2}{\Gamma \left( \alpha \right) }\left\vert \frac{d}{dt}\left[
\int_{0}^{t}f^{\prime }\left( r\right) \left( t-r\right) ^{-H+1/2}dr\right]
\right\vert \text{ by IBP} \\
&=&\frac{1}{\Gamma \left( \alpha \right) }\left\vert \left[
\int_{0}^{t}f^{\prime }\left( r\right) \left( t-r\right) ^{-H-1/2}dr\right]
\right\vert \\
&\leq &\frac{c_{H,q}\left\Vert f^{\prime }\right\Vert _{L^{q}}}{\Gamma
\left( \alpha \right) }<\infty \text{ by H\"{o}lder inequality and for }q>%
\frac{1}{1/2-H}.
\end{eqnarray*}%
In particular, $\left\Vert D_{0+}^{h+1/2}\left[ f\right] \right\Vert
_{L^{2}}\leq \left\Vert D_{0+}^{h+1/2}\left[ f\right] \right\Vert
_{L^{\infty }}\leq C_{H}\left\Vert f^{\prime }\right\Vert _{L^{q}}.$
\end{proof}

\section{Appendix 2}

We use a couple times the following formula: for $X_{i}$ centered Gaussian
random variable 
\begin{equation}
{\mathbb{E}}\left( \prod_{i=1}^{2n}X_{i}\right) =\sum_{\substack{ \left(
i_{1}^{1},i_{2}^{1}\right) ,...,\left( i_{1}^{n},i_{2}^{n}\right)  \\ \text{%
distinct pairs}}}\prod_{j=1}^{n}{\mathbb{E}}\left(
X_{i_{1}^{j}}X_{i_{2}^{j}}\right)  \label{Gaussian}
\end{equation}

\begin{proposition}
\label{increasing}For all $u<v$ in $\left[ s,t\right] ,$ $u^{\prime
}<v^{\prime }$ in $\left[ s^{\prime },t^{\prime }\right] ,$ where $s^{\prime
}<t^{\prime }\leq s<t$%
\begin{equation*}
0\leq -{\mathbb{E}}(W^{i}(u^{\prime },v^{\prime })W^{i}(u,v))\leq -{\mathbb{E%
}}(W^{i}(s^{\prime },t^{\prime })W^{i}(s,t))
\end{equation*}%
and for%
\begin{equation*}
0\leq {\mathbb{E}}(W^{i}(s,u)W^{i}(s,v))\leq {\mathbb{E}}%
(W^{i}(s,u)W^{i}(s,t))\leq {\mathbb{E}}\left( W^{i}(s,t)^{2}\right) .
\end{equation*}
\end{proposition}

\begin{proof}
The first inequality is obvious from the formula%
\begin{align*}
{\mathbb{E}}(W(u^{\prime },v^{\prime })W(u,v))& =\frac{1}{2}[|v-u^{\prime
}|^{2H}+|u-v^{\prime }|^{2H}-|u-u^{\prime }|^{2H}-|v-v^{\prime }|^{2H}] \\
& =H(2H-1)\int_{u}^{v}\int_{u^{\prime }}^{v^{\prime }}|y-x|^{2H-2}dydx.
\end{align*}%
We prove the second inequality with basic inequalities:%
\begin{eqnarray*}
&&{\mathbb{E}}(W^{i}(s,u)W^{i}(s,v)) \\
&=&\frac{1}{2}\left[ |u-s|^{2H}+|v-s|^{2H}-|v-u|^{2H}\right] \\
&\leq &\frac{1}{2}\left[ \left\vert u-s\right\vert ^{2H}+\left\vert
t-s\right\vert ^{2H}-|t-u|^{2H}\right] ={\mathbb{E}}(W^{i}(s,u)W^{i}(s,t)) \\
&\leq &\frac{1}{2}\left[ \left\vert t-s\right\vert ^{2H}+\left\vert
t-s\right\vert ^{2H}-0\right] ={\mathbb{E}}\left( W^{i}(s,t)^{2}\right) .
\end{eqnarray*}
\end{proof}

In the lemmae below, we will use the notation of the proof of lemma \ref%
{leadtothemainth}.

\begin{lemma}
\label{1lem}Let $f_{u}$ and $g_{u}$ two measurable processes independent of $%
B,$ and Skorokhod integrable. Then, for all $s^{\prime }<t^{\prime }\leq s<t$%
,%
\begin{equation*}
{\mathbb{E}}\left( \left( \int_{s}^{t}f\left( u\right) dB_{u}\right) \left(
\int_{s^{\prime }}^{t^{\prime }}g\left( u\right) dB_{u}\right) \right) =0.
\end{equation*}
\end{lemma}

\begin{proof}
Obvious from the property of Skorokhod integrals.
\end{proof}

\begin{lemma}
\label{2lem}For all $s^{\prime }<t^{\prime }\leq s<t$,%
\begin{equation*}
{\mathbb{E}}\left( \left( \mathbf{A}_{s^{\prime },t^{\prime }}^{2}\right)
^{i,j}\left( \mathbf{A}_{s,t}^{2}\right) ^{i,j}\right) \leq {\mathbb{E}}%
(W^{i}(s,t)W^{i}(s^{\prime },t^{\prime }))^{2}.
\end{equation*}
\end{lemma}

\begin{proof}
Once again, we can take $i\neq j$. Hence, denoting $X=\left\vert {\mathbb{E}}%
\left( \left( \mathbf{A}_{s^{\prime },t^{\prime }}^{2}\right) ^{i,j}\left( 
\mathbf{A}_{s,t}^{2}\right) ^{i,j}\right) \right\vert ,$ we have that $X$
equals%
\begin{gather*}
\left\vert {\mathbb{E}}\left( \int_{0}^{s}\left( \int_{s}^{t}\partial
K_{H}(r,u)W_{H}^{i}(s,r)dr\right) dB_{u}^{j}\int_{0}^{s^{\prime }}\left(
\int_{s^{\prime }}^{t^{\prime }}\partial K_{H}(r^{\prime
},u)W_{H}^{i}(s^{\prime },r^{\prime })dr^{\prime }\right) dB_{u}^{j}\right)
\right\vert \\
=\int_{u=0}^{s^{\prime }}\int_{r=s}^{t}\int_{r^{\prime }=s^{\prime
}}^{t^{\prime }}\partial K_{H}(r,u)\partial K_{H}(r^{\prime },u){\mathbb{E}}%
\left( W^{i}(s,r)W^{i}(s^{\prime },r^{\prime })\right) drdr^{\prime }du.
\end{gather*}%
Now we bound ${\mathbb{E}}\left( W_{H}^{i}(s,r)W_{H}^{i}(s^{\prime
},r^{\prime })\right) $ by $|{\mathbb{E}}\left( W^{i}(s,t)W^{i}(s^{\prime
},t^{\prime })\right) |$, thanks to proposition \ref{increasing}. Hence, we
obtain that%
\begin{eqnarray*}
X &\leq &\left\vert {\mathbb{E}}\left( W^{i}(s,t)W^{i}(s^{\prime },t^{\prime
})\right) \right\vert \int_{u=0}^{s^{\prime }}\int_{r=s}^{t}\int_{r^{\prime
}=s^{\prime }}^{t^{\prime }}\partial K_{H}(r,u)\partial K_{H}(r^{\prime
},u)drdr^{\prime }du \\
&=&\left\vert {\mathbb{E}}\left( W^{i}(s,t)W^{i}(s^{\prime },t^{\prime
})\right) \right\vert \int_{u=0}^{s^{\prime }}\left[ K_{H}(t,u)-K(s,u)\right]
\left[ K_{H}(t^{\prime },u)-K(s^{\prime },u)\right] du.
\end{eqnarray*}

Then, 
\begin{eqnarray*}
\int_{u=0}^{s^{\prime }}\left[ K_{H}(t,u)-K(s,u)\right] K_{H}(t^{\prime
},u)du &\leq &\int_{u=0}^{t^{\prime }}\left[ K_{H}(t,u)-K(s,u)\right]
K_{H}(t^{\prime },u)du \\
&=&{\mathbb{E}}\left( W^{i}(s,t)W^{i}(t^{\prime })\right)
\end{eqnarray*}%
and%
\begin{equation*}
\int_{u=0}^{s^{\prime }}\left[ K_{H}(t,u)-K(s,u)\right] K(s^{\prime },u)du={%
\mathbb{E}}\left( W^{i}(s,t)W^{i}(s^{\prime })\right) .
\end{equation*}%
That gives our result.
\end{proof}

The same argument gives the following lemma:

\begin{lemma}
\label{2lemcopy}For all $s<t$,%
\begin{equation*}
{\mathbb{E}}\left( \left\vert \left( \mathbf{A}_{s,t}^{2}\right)
^{i,j}\right\vert ^{2}\right) \leq {\mathbb{E}}(\left\vert
W^{i}(s,t)\right\vert ^{2})^{2},
\end{equation*}
\end{lemma}

\begin{lemma}
\label{3lem}For all $s^{\prime }<t^{\prime }\leq s<t$,%
\begin{equation*}
{\mathbb{E}}\left( \left( \mathbf{Y}_{s^{\prime },t^{\prime }}\right)
^{i,j}\left( \mathbf{Y}_{s,t}\right) ^{i,j}\right) \leq {\mathbb{E}}%
(W^{i}(s,t)W^{i}(s^{\prime },t^{\prime }))^{2}
\end{equation*}
\end{lemma}

\begin{proof}
Recall that 
\begin{equation*}
\mathbf{Y}_{s,t}=\frac{1}{|t-s|^{2}}\int_{s}^{t}[(u-s)W(u,t)-(t-u)W(s,u)]du%
\otimes W(s,t)
\end{equation*}

Moreover, since 
\begin{equation*}
{\mathbb{E}}(W(s)W(t))=\frac{1}{2}[s^{2H}+t^{2H}-|t-s|^{2H}]
\end{equation*}%
we obtain ${\mathbb{E}}(\mathbf{Y}_{s,t})=0.$ Therefore, $%
2|t-s|^{2}|t^{\prime }-s^{\prime }|^{2}{\mathbb{E}}\left( \left( \mathbf{Y}%
_{s,t}\right) ^{i,j}\left( \mathbf{Y}_{s^{\prime },t^{\prime }}\right)
^{i,j}\right) $ is bounded by $M_{1}M_{2}+M_{3}M_{4},$ where%
\begin{align*}
M_{1}& ={\mathbb{E}}\left(
\int_{s}^{t}[(u-s)W^{i}(u,t)-(t-u)W^{i}(s,u)]duW_{H}^{j}(s^{\prime
},t^{\prime })\right) , \\
M_{2}& ={\mathbb{E}}\left( \int_{s^{\prime }}^{t^{\prime }}[(u-s^{\prime
})W^{i}(u,t^{\prime })-(t^{\prime }-u)W^{i}(s^{\prime
},u)]duW_{H}^{j}(s,t)\right) , \\
M_{3}& ={\mathbb{E}}\left( \int_{s}^{t}[(u-s)W^{i}(u,t)-(t-u)W^{i}(s^{\prime
},u)]du\right. \\
& \left. \times \int_{s^{\prime }}^{t^{\prime }}[(v-s^{\prime
})W^{j}(u,t^{\prime })-(t^{\prime }-v)W^{j}(s^{\prime },v)]dv\right) , \\
M_{4}& ={\mathbb{E}}\left( W^{j}(s,t)W^{j}(s^{\prime },t^{\prime })\right) ,
\end{align*}

Using proposition \ref{increasing}, we obtain that%
\begin{eqnarray*}
M_{1} &\leq &\int_{s}^{t}\left( (u-s)+\left( t-u\right) \right) du\left\vert 
{\mathbb{E}}\left( W^{i}(s,t)W^{j}(s^{\prime },t^{\prime })\right)
\right\vert \\
&\leq &\left( t-s\right) ^{2}\left\vert {\mathbb{E}}\left(
W^{i}(s,t)W^{j}(s^{\prime },t^{\prime })\right) \right\vert .
\end{eqnarray*}%
Similarly, we have%
\begin{eqnarray*}
M_{2} &\leq &\left( t^{\prime }-s^{\prime }\right) ^{2}\left\vert {\mathbb{E}%
}\left( W^{i}(s,t)W^{j}(s^{\prime },t^{\prime })\right) \right\vert , \\
M_{3} &\leq &\left( t-s\right) ^{2}\left( t^{\prime }-s^{\prime }\right)
^{2}\left\vert {\mathbb{E}}\left( W^{i}(s,t)W^{i}(s^{\prime },t^{\prime
})\right) \right\vert .
\end{eqnarray*}%
Hence, we have shown that%
\begin{equation*}
\left\vert {\mathbb{E}}\left( \left( \mathbf{Y}_{s,t}\right) ^{i,j}\left( 
\mathbf{Y}_{s^{\prime },t^{\prime }}\right) ^{i,j}\right) \right\vert \leq {%
\mathbb{E}}\left( W^{i}(s,t)W^{i}(s^{\prime },t^{\prime })\right) ^{2}
\end{equation*}
\end{proof}

\begin{lemma}
\label{3lembis}For all $s^{\prime }<t^{\prime }\leq s<t$ or $\left(
s^{\prime },t^{\prime }\right) =\left( s,t\right) ,$%
\begin{eqnarray*}
&&{\mathbb{E}}\left( \left( \mathbf{\tilde{Y}}_{s^{\prime },t^{\prime
}}\right) ^{i,j}\left( \mathbf{\tilde{Y}}_{s,t}\right) ^{i,j}\right) \\
&\leq &C{\mathbb{E}}\left( W^{j}(s,t)W^{j}(s^{\prime },t^{\prime })\right) {%
\mathbb{E}}\left( \left\vert W^{j}(s,t)\right\vert ^{2}\right) {\mathbb{E}}%
\left( \left\vert W^{j}(s^{\prime },t^{\prime })\right\vert ^{2}\right)
\end{eqnarray*}
\end{lemma}

\begin{proof}
Recall that 
\begin{equation*}
\mathbf{\tilde{Y}}_{s,t}=\frac{1}{|t-s|^{3}}%
\int_{s<u_{1}<u_{2}<t}[(u_{1}-s)W(u_{1},t)-(t-u_{1})W(s,u_{1})]du_{1}du_{2}%
\otimes W(s,t)^{\otimes 2}.
\end{equation*}%
The proof is quite similar to the previous one. By equality (\ref{Gaussian}),%
\begin{eqnarray*}
\left( t-s\right) ^{3}\left( t^{\prime }-s^{\prime }\right) ^{3}E\left( 
\mathbf{\tilde{Y}}_{s,t}\mathbf{\tilde{Y}}_{s^{\prime },t^{\prime }}\right)
&=&\tilde{M}_{3}\left( \tilde{M}_{5}\tilde{M}_{6}+2\tilde{M}_{4}^{2}\right)
\\
&&+2\tilde{M}_{7}\left( 2\tilde{M}_{8}\tilde{M}_{4}+\tilde{M}_{2}\tilde{M}%
_{6}\right) \\
&&+2\tilde{M}_{1}\left( 2\tilde{M}_{2}\tilde{M}_{4}+\tilde{M}_{8}\tilde{M}%
_{5}\right)
\end{eqnarray*}%
where%
\begin{align*}
\tilde{M}_{1}& ={\mathbb{E}}\left( \left(
\int_{s<u_{1}<u_{2}<t}[(u_{1}-s)W(u_{1},t)-(t-u_{1})W(s,u_{1})]du_{1}du_{2}%
\right) W^{j}(s^{\prime },t^{\prime })\right) , \\
\tilde{M}_{2}& ={\mathbb{E}}\left( \left( \int_{s^{\prime
}<u_{1}<u_{2}<t^{\prime }}[(u_{1}-s^{\prime })W^{i}(u_{1},t^{\prime
})-(t^{\prime }-u_{1})W^{i}(s^{\prime },u_{1})]du_{1}du_{2}\right)
W^{j}(s,t)\right) , \\
\tilde{M}_{3}& ={\mathbb{E}}\left(
\int_{s<u_{1}<u_{2}<t}[(u_{1}-s)W^{i}(u_{1},t)-(t-u_{1})W^{i}(u_{1},t)]du_{1}du_{2}\right.
\\
& \left. \int_{s^{\prime }<u_{1}<u_{2}<t^{\prime }}[(u_{1}-s^{\prime
})W^{i}(u_{1},t^{\prime })-(t^{\prime }-u_{1})W^{i}(u_{1},t^{\prime
})]du_{1}du_{2}\right) , \\
\tilde{M}_{4}& ={\mathbb{E}}\left( W^{j}(s,t)W^{j}(s^{\prime },t^{\prime
})\right) , \\
\tilde{M}_{5}& ={\mathbb{E}}\left( \left\vert W^{j}(s,t)\right\vert
^{2}\right) , \\
\tilde{M}_{6}& ={\mathbb{E}}\left( \left\vert W^{j}(s^{\prime },t^{\prime
})\right\vert ^{2}\right) , \\
\tilde{M}_{7}& ={\mathbb{E}}\left( \left(
\int_{s<u_{1}<u_{2}<t}[(u_{1}-s)W(u_{1},t)-(t-u_{1})W(s,u_{1})]du_{1}du_{2}%
\right) W^{j}(s,t)\right) , \\
\tilde{M}_{8}& ={\mathbb{E}}\left( \left( \int_{s^{\prime
}<u_{1}<u_{2}<t^{\prime }}[(u_{1}-s^{\prime })W^{i}(u_{1},t^{\prime
})-(t^{\prime }-u_{1})W^{i}(s^{\prime },u_{1})]du_{1}du_{2}\right)
W^{j}(s^{\prime },t^{\prime })\right) .
\end{align*}%
Working just as in the previous lemma , we see that%
\begin{equation*}
\frac{\left\vert \tilde{M}_{1}\right\vert }{\left( t-s\right) ^{3}}+\frac{%
\left\vert \tilde{M}_{2}\right\vert }{\left( t^{\prime }-s^{\prime }\right)
^{3}}+\frac{\left\vert \tilde{M}_{3}\right\vert }{\left( t-s\right)
^{3}\left( t^{\prime }-s^{\prime }\right) ^{3}}\leq C\tilde{M}_{4},
\end{equation*}%
and $\left\vert \tilde{M}_{7}\right\vert \leq C\left( t-s\right) ^{3}\tilde{M%
}_{5}$ and $\left\vert \tilde{M}_{8}\right\vert \leq C\left( t-s\right) ^{3}%
\tilde{M}_{6}$. Therefore, we obtain that%
\begin{equation*}
{\mathbb{E}}\left( \left( \mathbf{\tilde{Y}}_{s^{\prime },t^{\prime
}}\right) ^{i,j}\left( \mathbf{\tilde{Y}}_{s,t}\right) ^{i,j}\right) \leq
C\left( \tilde{M}_{4}^{3}+\tilde{M}_{4}\tilde{M}_{5}\tilde{M}_{6}\right) .
\end{equation*}%
Finally, by Cauchy-Schwarz, $\tilde{M}_{4}^{2}\leq \tilde{M}_{5}\tilde{M}%
_{6},$ which implies that%
\begin{equation*}
{\mathbb{E}}\left( \left( \mathbf{\tilde{Y}}_{s^{\prime },t^{\prime
}}\right) ^{i,j}\left( \mathbf{\tilde{Y}}_{s,t}\right) ^{i,j}\right) \leq C%
\tilde{M}_{4}\tilde{M}_{5}\tilde{M}_{6}.
\end{equation*}
\end{proof}

\begin{lemma}
\label{forthirdlevel}For $i\neq j,$ all $s^{\prime }<t^{\prime }\leq s<t$ or 
$s=s^{\prime }<t=t^{\prime },$ with $\left( s^{\prime },t^{\prime
},s,t\right) \in D,$%
\begin{multline*}
{\mathbb{E}}\left( \left( \frac{1}{t-s}\int_{s}^{t}\left( \mathbf{A}%
_{s,u}^{2}\right) ^{i,j}du\right) W^{i}\left( s,t\right) \left( \frac{1}{%
t^{\prime }-s^{\prime }}\int_{s^{\prime }}^{t^{\prime }}\left( \mathbf{A}%
_{s^{\prime },u}^{2}\right) ^{i,j}du\right) W^{i}\left( s^{\prime
},t^{\prime }\right) \right) \\
\leq C{\mathbb{E}}\left( W^{1}(s,t)W^{1}(s^{\prime },t^{\prime })\right) {%
\mathbb{E}}\left( \left\vert W^{1}(s,t)\right\vert ^{2}\right) {\mathbb{E}}%
\left( \left\vert W^{1}(s^{\prime },t^{\prime })\right\vert ^{2}\right) .
\end{multline*}
\end{lemma}

\begin{proof}
First observe that%
\begin{multline*}
\int_{s}^{t}\left( \mathbf{A}_{s,u}^{2}\right) ^{i,j}dW^{D,i}\left( u\right)
=\frac{W_{s,t}^{i}}{t-s}\int_{s}^{t}\left( \int_{0}^{s}\left(
\int_{s}^{u}\partial K_{H}(r,v)W^{i}(s,r)dr\right) dB_{v}^{j}\right) du \\
=\frac{1}{t-s}\int_{s}^{t}\left( \int_{0}^{s}\left(
W_{s,t}^{i}\int_{s}^{u}\partial K_{H}(r,v)W^{i}(s,r)dr\right)
dB_{v}^{j}\right) du.
\end{multline*}%
Therefore,%
\begin{align*}
& \left( t-s\right) \left( t^{\prime }-s^{\prime }\right) {\mathbb{E}}\left(
\int_{s}^{t}\left( \mathbf{A}_{s,u}^{2}\right) ^{i,j}dW^{D,i}\left( u\right)
\int_{s^{\prime }}^{t^{\prime }}\left( \mathbf{A}_{s,u}^{2}\right)
^{i,j}dW^{D,i}\left( u\right) \right) \\
& =\int_{s}^{t}\int_{s^{\prime }}^{t^{\prime }}{\mathbb{E}}\left( \left(
\int_{0}^{s}\left( W_{s,t}^{i}\int_{s}^{u}\partial
K_{H}(r,v)W_{s,r}^{i}dr\right) dB_{v}^{j}\right) \right. \\
& \,\,\,\,\,\,\,\,\left. \left( \int_{0}^{s^{\prime }}\left( W_{s^{\prime
},t^{\prime }}^{i}\int_{s^{\prime }}^{u^{\prime }}\partial K_{H}(r^{\prime
},v^{\prime })W_{s^{\prime },r^{\prime }}^{i}dr^{\prime }\right)
dB_{v^{\prime }}^{j}\right) \right) dudu^{\prime } \\
& =\int_{s}^{t}\int_{s^{\prime }}^{t^{\prime }}{\mathbb{E}}\left(
\int_{0}^{s^{\prime }}W_{s^{\prime },t^{\prime }}^{i}W_{s,t}^{i}\left(
\int_{s}^{u}\partial K_{H}(r,v)W_{s,r}^{i}dr\right) \right. \\
& \,\,\,\,\,\,\,\left. \left( \int_{s^{\prime }}^{u^{\prime }}\partial
K_{H}(r^{\prime },v)W_{s^{\prime },r^{\prime }}^{i}dr^{\prime }\right)
dv\right) dudu^{\prime } \\
& =\int_{s}^{t}\int_{s^{\prime }}^{t^{\prime }}\int_{0}^{s^{\prime
}}\int_{s}^{u}\int_{s^{\prime }}^{u^{\prime }}\partial K_{H}(r^{\prime
},v)\partial K_{H}(r,v){\mathbb{E}}\left( W_{s^{\prime },t^{\prime
}}^{i}W_{s,t}^{i}W_{s,r}^{i}W_{s^{\prime },r^{\prime }}^{i}\right)
drdr^{\prime }dvdudu^{\prime }.
\end{align*}%
Now, using equality (\ref{Gaussian}) and proposition \ref{increasing}, we
see that 
\begin{eqnarray*}
{\mathbb{E}}\left( W_{s^{\prime },t^{\prime
}}^{i}W_{s,t}^{i}W_{s,r}^{i}W_{s^{\prime },r^{\prime }}^{i}\right) &\leq &2{%
\mathbb{E}}\left( W_{s^{\prime },t^{\prime }}^{i}W_{s,t}^{i}\right) ^{2}+{%
\mathbb{E}}\left( \left\vert W_{s,t}^{i}\right\vert ^{2}\right) {\mathbb{E}}%
\left( \left\vert W_{s^{\prime },t^{\prime }}^{i}\right\vert ^{2}\right) \\
&\leq &3{\mathbb{E}}\left( \left\vert W_{s,t}^{i}\right\vert ^{2}\right) {%
\mathbb{E}}\left( \left\vert W_{s^{\prime },t^{\prime }}^{i}\right\vert
^{2}\right) \text{ by Cauchy-Schwartz.}
\end{eqnarray*}%
Hence%
\begin{eqnarray*}
&&\left( t-s\right) \left( t^{\prime }-s^{\prime }\right) {\mathbb{E}}\left(
\int_{s}^{t}\left( \mathbf{A}_{s,u}^{2}\right) ^{i,j}dW^{D,i}\left( u\right)
\int_{s^{\prime }}^{t^{\prime }}\left( \mathbf{A}_{s,u}^{2}\right)
^{i,j}dW^{D,i}\left( u\right) \right) \\
&\leq &3{\mathbb{E}}\left( \left\vert W_{s,t}^{i}\right\vert ^{2}\right) {%
\mathbb{E}}\left( \left\vert W_{s^{\prime },t^{\prime }}^{i}\right\vert
^{2}\right) \\
&&\times \int_{u=s}^{t}\int_{u^{\prime }=s^{\prime }}^{t^{\prime
}}\int_{v=0}^{s^{\prime }}\int_{r=s}^{u}\int_{r^{\prime }=s^{\prime
}}^{u^{\prime }}\partial K_{H}(r^{\prime },v)\partial K_{H}(r,v)drdr^{\prime
}dvdudu^{\prime }.
\end{eqnarray*}%
By positivity of the derivative of the kernel $K_{H}$, we therefore obtain%
\begin{eqnarray*}
&&\left( t-s\right) \left( t^{\prime }-s^{\prime }\right) {\mathbb{E}}\left(
\int_{s}^{t}\left( \mathbf{A}_{s,u}^{2}\right) ^{i,j}dW^{D,i}\left( u\right)
\int_{s^{\prime }}^{t^{\prime }}\left( \mathbf{A}_{s,u}^{2}\right)
^{i,j}dW^{D,i}\left( u\right) \right) \\
&\leq &3{\mathbb{E}}\left( \left\vert W_{s,t}^{i}\right\vert ^{2}\right) {%
\mathbb{E}}\left( \left\vert W_{s^{\prime },t^{\prime }}^{i}\right\vert
^{2}\right) \left( t-s\right) \left( t^{\prime }-s^{\prime }\right) \\
&&\times \int_{v=0}^{s^{\prime }}\int_{r=s}^{t}\int_{r^{\prime }=s^{\prime
}}^{t^{\prime }}\partial K_{H}(r^{\prime },v)\partial K_{H}(r,v)drdr^{\prime
}dv.
\end{eqnarray*}%
We proved in lemma \ref{2lem} that 
\begin{equation*}
\int_{v=0}^{s^{\prime }}\int_{r=s}^{t}\int_{r^{\prime }=s^{\prime
}}^{t^{\prime }}\partial K_{H}(r^{\prime },v)\partial K_{H}(r,v)drdr^{\prime
}dv\leq {\mathbb{E}}\left( W^{i}(s,t)W^{i}(s^{\prime },t^{\prime })\right) ,
\end{equation*}
which conclude the case $s^{\prime }<t^{\prime }\leq s<t.$ We leave the
easier case $s=s^{\prime }<t=t^{\prime }$ to the reader.
\end{proof}

\begin{lemma}
\label{2lemdegree3}For all $x\leq s^{\prime }<t^{\prime }\leq s<t$, and $%
k\notin \left\{ i,j\right\} $%
\begin{eqnarray*}
&&{\mathbb{E}}\left( W_{x,s^{\prime }}^{i}\left( \mathbf{A}_{s^{\prime
},t^{\prime }}^{2}\right) ^{j,k}W_{x,s}^{i}\left( \mathbf{A}%
_{s,t}^{2}\right) ^{j,k}\right) \\
&\leq &{\mathbb{E}}(W_{x,s}^{1}W_{x,s^{\prime
}}^{1}W^{1}(s,t)W^{1}(s^{\prime },t^{\prime })){\mathbb{E}}\left(
W^{1}(s,t)W^{1}(s^{\prime },t^{\prime })\right) .
\end{eqnarray*}%
and%
\begin{equation*}
{\mathbb{E}}\left( \left\vert W_{x,s}^{i}\left( \mathbf{A}_{s,t}^{2}\right)
^{j,k}\right\vert ^{2}\right) \leq {\mathbb{E}}(\left\vert
W_{x,y}^{i}\right\vert ^{2}\left\vert W^{i}(s,t)\right\vert ^{2}){\mathbb{E}}%
\left( \left\vert W^{i}(s,t)\right\vert ^{2}\right)
\end{equation*}
\end{lemma}

\begin{proof}
Using the notation $X=\left\vert {\mathbb{E}}\left( W_{x,s^{\prime
}}^{i}\left( \mathbf{A}_{s^{\prime },t^{\prime }}^{2}\right)
^{j,k}W_{x,s}^{i}\left( \mathbf{A}_{s,t}^{2}\right) ^{j,k}\right)
\right\vert ,$ we have that $X$ equals%
\begin{gather*}
\begin{array}{l}
\left\vert {\mathbb{E}}\left( \int_{0}^{s}\left( \int_{s}^{t}\partial
K_{H}(r,u)W_{x,s}^{i}W^{j}(s,r)dr\right) dB_{u}^{j}\right. \right. \\ 
\text{ \ \ \ \ \ \ \ \ \ \ \ \ \ \ \ \ \ \ \ \ \ \ \ \ }\left. \left.
\int_{0}^{s^{\prime }}\left( \int_{s^{\prime }}^{t^{\prime }}\partial
K_{H}(r^{\prime },u)W_{x,s^{\prime }}^{i}W^{j}(s^{\prime },r^{\prime
})dr^{\prime }\right) dB_{u}^{j}\right) \right\vert%
\end{array}
\\
=\int_{u=0}^{s^{\prime }}\int_{r=s}^{t}\int_{r^{\prime }=s^{\prime
}}^{t^{\prime }}\partial K_{H}(r,u)\partial K_{H}(r^{\prime },u){\mathbb{E}}%
\left( W_{x,s}^{i}W_{x,s^{\prime }}^{i}W^{j}(s,r)W^{j}(s^{\prime },r^{\prime
})\right) drdr^{\prime }du.
\end{gather*}%
Now we bound ${\mathbb{E}}\left( W_{x,s}^{i}W_{x,s^{\prime
}}^{i}W^{j}(s,r)W^{j}(s^{\prime },r^{\prime })\right) $ by ${\mathbb{E}}%
\left( W_{x,s}^{1}W_{x,s^{\prime }}^{1}W^{1}(s,t)W^{1}(s^{\prime },t^{\prime
})\right) $, thanks to equation (\ref{Gaussian}) and proposition \ref%
{increasing}. Hence, we obtain that%
\begin{equation*}
X\leq {\mathbb{E}}\left( W_{x,s}^{1}W_{x,s^{\prime
}}^{1}W^{1}(s,t)W^{1}(s^{\prime },t^{\prime })\right) \int_{u=0}^{s^{\prime
}}\int_{r=s}^{t}\int_{r^{\prime }=s^{\prime }}^{t^{\prime }}\partial
K_{H}(r,u)\partial K_{H}(r^{\prime },u)drdr^{\prime }du.
\end{equation*}%
We conclude just as in lemma \ref{2lem}. The second part of the lemma follow
in a similar way.
\end{proof}

\begin{lemma}
\label{finallemma}Let $D=\left( t_{k}\right) _{k\in \left\{ m,\cdots
,n\right\} }$ be a subdivision. Then, for $p^{\prime }>1/H,$, there exists a
constant $C_{p^{\prime }}$ such that 
\begin{eqnarray*}
&&\sum_{k,l=m}^{n}\left\Vert W^{j}(t_{k},t_{k+1})\right\Vert
_{L^{2}}\left\Vert W^{j}(t_{l},t_{l+1})\right\Vert _{L^{2}}{\mathbb{E}}%
\left( W^{i}(t_{k},t_{k+1})W^{i}(t_{l},t_{l+1})\right) \\
&\leq &C_{p^{\prime }}\left( t_{n}-t_{m}\right) ^{4/p^{\prime
}}|D|^{4H-4/p^{\prime }}.
\end{eqnarray*}%
In particular,%
\begin{equation*}
\sum_{k,l=m}^{n}{\mathbb{E}}\left(
W^{i}(t_{k},t_{k+1})W^{i}(t_{l},t_{l+1})\right) ^{2}\leq C_{p^{\prime
}}\left( t_{n}-t_{m}\right) ^{4/p^{\prime }}|D|^{4h-4/p^{\prime }}.
\end{equation*}
\end{lemma}

\begin{proof}
We observed in the proof of Proposition \ref{increasing} that, for $%
s^{\prime }<t^{\prime }\leq s<t,$%
\begin{equation*}
{\mathbb{E}}(W(s^{\prime },t^{\prime
})W(s,t))=C_{H}\int_{u=s}^{t}\int_{v=s^{\prime }}^{t^{\prime
}}(u-v)^{2H-2}dudv
\end{equation*}%
with $C_{H}=H(2H-1).$We bound the term $(u-v)^{2H-2}$ in the integral by $%
(u-t^{\prime })^{-1+\varepsilon }(t^{\prime }-v)^{2H-\varepsilon -1}$ (for a
fixed $\varepsilon \in \left( 0,2H\right) $ ) to obtain%
\begin{equation*}
|{\mathbb{E}}(W(s^{\prime },t^{\prime })W(s,t))|\leq
C_{H}\int_{u=s}^{t}\int_{v=s^{\prime }}^{t^{\prime }}(u-t^{\prime
})^{-1+\varepsilon }du(t^{\prime }-v)^{2H-\varepsilon -1}dv.
\end{equation*}%
Integrating with respect to $v\in \left[ s^{\prime },t^{\prime }\right] ,$
we obtain 
\begin{equation}
{\mathbb{E}}(W(s^{\prime },t^{\prime })W(s,t))\leq C_{H,\varepsilon
}\int_{s}^{t}(u-t^{\prime })^{-1+\varepsilon }du|t^{\prime }-s^{\prime
}|^{2H-\varepsilon }.  \label{1}
\end{equation}%
Multiplying inequality (\ref{1}) by $\left( t-s\right) ^{h}\left( t^{\prime
}-s^{\prime }\right) ^{h}$, we get%
\begin{multline*}
\left( t-s\right) ^{H}\left( t^{\prime }-s^{\prime }\right) ^{H}{\mathbb{E}}%
(W(s^{\prime },t^{\prime })W(s,t)) \\
\leq C_{H,\varepsilon }\left( \int_{s}^{t}(u-t^{\prime })^{-1+\varepsilon
}du\right) (t^{\prime }-s^{\prime })^{3H-\varepsilon }|t-s|^{H}.
\end{multline*}%
Fix $l\in \left\{ m,\cdots ,n\right\} ,$ and apply the above inequality to $%
\left( s^{\prime },t^{\prime },s,t\right) =\left(
t_{l},t_{l+1},t_{k},t_{k+1}\right) :$%
\begin{gather*}
\sum_{k=l+1}^{n}{\mathbb{E}}\left( W(t_{l},t_{l+1})W(t_{k},t_{k+1})\right)
^{2} \\
\leq C_{H,\varepsilon }\sum_{k=l+1}^{n}\left(
\int_{t_{k}}^{t_{k+1}}(u-t_{l+1})^{-1+\varepsilon }du\right)
(t_{l+1}-t_{l})^{3H-\varepsilon }|t_{k+1}-t_{k}|^{H}.
\end{gather*}%
We then bound $|t_{k+1}-t_{k}|^{H}$ by $\left\vert D\right\vert ^{H}$: 
\begin{eqnarray*}
&&\sum_{k=l+1}^{n}\left( t_{l+1}-t_{l}\right) ^{H}\left(
t_{k+1}-t_{k}\right) ^{H}{\mathbb{E}}\left(
W(t_{l},t_{l+1})W(t_{k},t_{k+1})\right) ^{2} \\
&\leq &C_{H,\varepsilon }\sum_{k=l+1}^{n}\left(
\int_{t_{k}}^{t_{k+1}}(u-t_{l+1})^{-1+\varepsilon }du\right)
(t_{l+1}-t_{l})^{3H-\varepsilon }|D|^{H} \\
&=&C_{H,\varepsilon }\left(
\int_{t_{l+1}}^{t_{n}}(u-t_{l+1})^{-1+\varepsilon }du\right)
(t_{l+1}-t_{l})^{3H-\varepsilon }|D|^{H} \\
&\leq &C_{H,\varepsilon }\left( t_{n}-t_{l+1}\right) ^{\varepsilon
}(t_{l+1}-t_{l})|D|^{4H-1-\varepsilon }
\end{eqnarray*}%
And therefore, we obtain%
\begin{eqnarray*}
&&\sum_{l=m}^{n-2}\sum_{k=l+1}^{n-1}\left( t_{l+1}-t_{l}\right) ^{H}\left(
t_{k+1}-t_{k}\right) ^{H}{\mathbb{E}}\left(
W^{i}(t_{k},t_{k+1})W^{i}(t_{l},t_{l+1})\right) \\
&\leq &C_{H,\varepsilon }\sum_{l=m}^{n-2}\left( t_{n}-t_{l+1}\right)
^{\varepsilon }(t_{l+1}-t_{l})|D|^{4H-1-\varepsilon }.
\end{eqnarray*}%
Comparing $\sum_{l=m}^{n-2}\left( t_{n}-t_{l+1}\right) ^{\varepsilon
}(t_{l+1}-t_{l})$ to $\int_{t_{m}}^{t_{n}}\left( t_{n}-u\right)
^{\varepsilon }du$, we see that%
\begin{equation*}
\sum_{l=m}^{n-2}\left( t_{n}-t_{l+1}\right) ^{\varepsilon
}(t_{l+1}-t_{l})\leq C\left( t_{n}-t_{m}\right) ^{1+\varepsilon }.
\end{equation*}
Hence,%
\begin{eqnarray*}
&&\sum_{l=m}^{n-2}\sum_{k=l+1}^{n-1}\left( t_{l+1}-t_{l}\right) ^{H}\left(
t_{k+1}-t_{k}\right) ^{H}{\mathbb{E}}\left(
W^{i}(t_{k},t_{k+1})W^{i}(t_{l},t_{l+1})\right) \\
&\leq &C_{H,\varepsilon }\left( t_{n}-t_{m}\right) ^{1+\varepsilon
}|D|^{4H-1-\varepsilon }
\end{eqnarray*}%
Also, we easily obtain%
\begin{eqnarray*}
\sum_{l=m}^{n-1}\left( t_{l+1}-t_{l}\right) ^{2H}{\mathbb{E}}\left(
W(t_{l},t_{l+1})^{2}\right) &=&C_{H}\sum_{l=m}^{n-1}\left(
t_{l+1}-t_{l}\right) ^{4H} \\
&\leq &C_{H}|D|^{4H-1-\varepsilon }\sum_{l=m}^{n-1}\left(
t_{l+1}-t_{l}\right) ^{1+\varepsilon } \\
&\leq &C_{H,\varepsilon }\left( t_{n}-t_{m}\right) ^{1+\varepsilon
}|D|^{4H-1-\varepsilon }
\end{eqnarray*}%
Adding the two previous inequalities give our results:%
\begin{multline*}
\sum_{k,l=m}^{n}\left\Vert W^{j}(t_{k},t_{k+1})\right\Vert
_{L^{2}}\left\Vert W^{j}(t_{l},t_{l+1})\right\Vert _{L^{2}}{\mathbb{E}}%
\left( W^{i}(t_{k},t_{k+1})W^{i}(t_{l},t_{l+1})\right) \\
\leq C_{H,\varepsilon }\left( t_{n}-t_{m}\right) ^{1+\varepsilon
}|D|^{4H-1-\varepsilon }
\end{multline*}%
which is our result with $1+\varepsilon =4/p^{\prime }$ (which is coherent
with $\varepsilon \in \left( 0,2H\right) $, as $H\in \left( 1/4,1/2\right) $%
).
\end{proof}

\bigskip


\begin{thebibliography}{99}
\bibitem{F} Friz, P: Continuity of the It\^{o}-map for H\"{o}lder rough path
with applications to the Support Theorem in H\"{o}lder norm (To appear at
IMA Proceedings)
\end{thebibliography}
\end{document}